\documentclass[a4paper,10pt]{article}
\usepackage{amsmath}
\usepackage{amssymb,amsfonts,amsthm}
\usepackage[applemac]{inputenc}
\usepackage[english]{babel}
\usepackage{caption}
\usepackage{subcaption}
\usepackage{anysize}
\usepackage{enumerate}
\usepackage{tikz}
\usepackage[all]{xy}
\usepackage{mathrsfs}
\usepackage{verbatim}
\usepackage{cite}
\usepackage{functan} 
\newcommand{\abs}[1]{\left\vert#1\right\vert}

\usepackage[width=0.8\textwidth]{caption}

\renewenvironment{abstract}{\small\quotation\noindent
 {\bfseries \abstractname .}}{\endquotation \par}

\newenvironment{trivlistparm}[2]{\trivlist
\item[{#1 #2}]}{\endtrivlist}

\newenvironment{prooftext}[1]{\trivlistparm{\bfseries}{#1}}{\Qed\endtrivlistparm}
\newenvironment{demo}{\trivlistparm{\bfseries}{Proof.}}{\Qed\endtrivlistparm}

\newcounter{claimcounter}
\newenvironment{claim}{\stepcounter{claimcounter}{\par\noindent\underline{Claim \theclaimcounter:}}}{}

\catcode`\@=11

\def\resetthefootnote{\renewcommand{\thefootnote}{\@arabic\c@footnote} }
\def\@principiremex#1{\trivlist
 \item[\hskip \labelsep{\bfseries #1\ \thethm.}]\ignorespaces}
\def\opar@principiremex#1[#2]{\trivlist
 \item[\hskip \labelsep{\bfseries #1\ \thethm\ (#2).}]\ignorespaces}

\newcommand{\newTHEOremrom}[2]{\newenvironment{#1}{\refstepcounter{thm}\@ifnextchar[{\opar@principiremex{#2}}
{\@principiremex{#2}}}{\qedB\endtrivlist}} \catcode`\@=12
\DeclareMathSymbol{\square}{\mathord}{AMSa}{"03}
\newcommand{\qedB}{\nopagebreak\hspace*{\fill}$\square$\par}
\newcommand{\Qed}{\nopagebreak\hspace*{\fill}{\vrule width6pt height6pt depth0pt}\par}

\newTHEOremrom{defi} {Definition}
\newTHEOremrom {rem} {Remark}
\newTHEOremrom {ex} {Example}

\renewcommand{\geq}{\geqslant}
\renewcommand{\epsilon}{\varepsilon}
\renewcommand{\leq}{\leqslant}
\newcommand{\R}{\mathbb{R}}

\newcommand{\Z}{\mathbb{Z}}

\newcommand{\C}{\mathbb{C}}

\newcommand{\T}{\mathbb{T}}

\newtheorem{thm}{Theorem}[section]
\newtheorem{thmx}{Theorem}
\newtheorem{prop}[thm]{Proposition}
\newtheorem{cor}[thm]{Corollary}
\newtheorem{lema}[thm]{Lemma}

\title{\textbf{Periodic oscillators, isochronous centers and resonance}
\footnotetext{2010 {\it Mathematics Subject Classification.} 34C15; 34C10, 34D05, 34D10, 34D23.}
\footnotetext{{\it Key words and phrases}:  Isochronous center, oscillator, resonance, perturbation.}
\footnotetext{All the authors are partially supported by the MINECO/FEDER grant MTM2017-82348-C2-1-P. D. Rojas is also partially supported by the MINECO/FEDER grant MTM2017-86795-C3-1-P.}
\footnotetext{{\it Email addresses:} \texttt{rojas@ugr.es}  (D.~Rojas, corresponding author), \texttt{rortega@ugr.es} (R.~Ortega).}
}

\author{Rafael Ortega and David Rojas \\[10pt]
{\small \textsl{Departamento de Matem\'atica Aplicada,}}\\
\vspace{-2pt}
{\small \textsl{Universidad de Granada, 18071 Granada, Spain}}}

\date{}

\begin{document}

\maketitle

\begin{abstract}
An oscillator is called isochronous if all motions have a common period. When the system is forced by a time-dependent perturbation with the same period the dynamics may change and the phenomenon of resonance can appear. In this context, resonance means that all solutions are unbounded. The theory of resonance is well known for the harmonic oscillator and we extend it to nonlinear isochronous oscillators.
\end{abstract}

\section{Introduction}\label{sec:intro}
Consider an oscillator with equation
\begin{equation}\label{aut}
\ddot x + V'(x)=0,\ x\in\R
\end{equation}
and assume that it has an isochronous center at the origin. This means that $x=0$ is the only equilibrium of the equation and the remaining solutions are periodic with a fixed period, say $T=2\pi$. We are interested in the phenomenon of resonance for periodic perturbations. More precisely, we ask for the class of $2\pi$-periodic functions $p(t)$ such that all the solutions of the non-autonomous equation
\begin{equation}\label{per}
\ddot x + V'(x)=\epsilon p(t)
\end{equation}
are unbounded. Here $\epsilon\neq 0$ is a small parameter.

The simplest isochronous center is produced by the harmonic oscillator, $V(x)=\frac{1}{2}n^2x^2$, $n=1,2,\dots$ In this case the previous question has a well-known answer: resonance occurs whenever the integral
\[
I_n(p)\!:=\int_0^{2\pi} p(t)e^{int}dt
\]
does not vanish. After this example the study of resonance for general isochronous oscillators seems natural. As far as we know this question was first raised by Prof. Roussarie in the Open Problems Session of the II Symposium on Planar Vector Fields (Lleida, 2000). Concrete examples of functions $p(t)$ producing resonance were presented in \cite{Ortega}. See also \cite{BFS}. 

The goal of the present paper is to identify a general class of forcings leading to resonance. Our main result can be interpreted as a nonlinear version of the condition $I_n(p)\neq 0$. To explain this in more precise terms we consider a cylinder $\mathcal C=(\R/2\pi\Z)\times[0,\infty)$ with coordinates $(\theta,r)$. The solution of \eqref{aut} with initial conditions $x(0)=r$, $\dot x(0)=0$ will be denoted by $\varphi(t,r)$. The complex-valued solution of the linear problem
\[
\ddot y+ V''(\varphi(t,r))y = 0,\ y(0)=1,\ \dot y(0)=i
\]
is denoted by $\psi(t,r)$. Then we define a function on the cylinder
\[
\Phi_p:\mathcal C\rightarrow\C,\ \Phi_p(\theta,r)=\frac{1}{2\pi}\int_0^{2\pi} p(t-\theta)\psi(t,r)dt
\]
and prove, under certain assumptions on the potential $V$, that all solutions of \eqref{per} are unbounded when the condition
\begin{equation}\label{rc}
\inf_{\mathcal C}\abs{\Phi_p(\theta,r)}>0
\end{equation}
holds.

This is a sufficient condition for resonance but it is not too far from being also necessary. A partial converse of the main theorem holds: a periodic solution exists when the function $\Phi_p$ has a non-degenerate zero. Note that periodic solutions are bounded and so resonance is excluded.

Let us now discuss the form of $\Phi_p$ in some particular cases. For the linear oscillator, $V(x)=\frac{1}{2}n^2x^2$, $n=1,2,\dots$ elementary computations lead to the estimates
\[
\frac{1}{2\pi n} |I_n (p)|\leq |\Phi_p (\theta ,r)|\leq \frac{1}{2\pi} |I_n (p)|.
\]

The condition $I_n(p)\neq 0$ is equivalent to \eqref{rc}. Another example that has been widely studied is the asymmetric oscillator associated to the potential
\[
V(x)=\frac{1}{2}\left(\alpha(x^+)^2+\beta(x^-)^2\right)
\]
where $x^+=\max\{x,0\}$, $x^-=\max\{-x,0\}$ and $\alpha$, $\beta$ are positive constants with $\frac{1}{\sqrt{\alpha}}+\frac{1}{\sqrt{\beta}}=\frac{2\pi}{n}$, $n=1,2,\dots$ In this case the homogeneity of the equation \eqref{aut} implies that $\Phi_p(\theta,r)=\Phi_p(\theta,1)$. The function $\Phi_p(\theta,1)$ appeared first in the work of Dancer on the periodic problem (see~\cite{Dancer}). Assuming that the zeros of the function $\Phi_p(\theta,1)$ were non-degenerate, it was proved in~\cite{AO} that all solutions with large initial condition were unbounded. Soon after, it was proved in~\cite{Liu} that all solution are bounded when $\Phi_p(\theta,1)$ does not vanish. This boundedness condition improved a previous result in~\cite{Ortega96} of more local nature. The boundedness of all solutions has been also treated in~\cite{BonFab,LiJin}. These papers deal with a class of isochronous oscillators with the same asymptotic behavior as the asymmetric oscillator.

For a general isochronous oscillator checking the condition \eqref{rc} can be difficult. We have analysed in detail the equation
\[
\ddot x + \frac{1}{4}\left( x+1 - \frac{1}{(x+1)^3}\right) = \epsilon p(t), \ x\in(1,+\infty).
\]
In this case the functions $\varphi(t,r)$ and $\psi(t,r)$ can be obtained explicitly and the resonance condition~\eqref{rc} can be reformulated in terms of the Fourier coefficients of $p(t)$. In particular we obtain an extension of the result in~\cite{BFS} for $p(t)=\sin t$. This equation has a singularity at $x=-1$ and the motion is constrained to a proper and unbounded interval. It can be seen as a prototype of the class of isochronous oscillators with one asymptote. Incidentally we note that there are no isochronous oscillators having two asymptotes so that the motion is constrained to a bounded interval. See \cite{OrtRoj} for more details. Some boundedness results for isochronous oscillators with singularities can be found in \cite{CapDamLiu, NewLiu}.

Up to now all non-autonomous perturbations have been of additive type. This is very natural if we have in mind mechanical oscillators but different perturbations can appear in other contexts. Motivated by a geometrical problem, Ai, Chou and Wei considered in \cite{ACW} the equation
\begin{equation}\label{acw}
\ddot x+x=\frac{R(t)}{x^3},\ x>0
\end{equation}
where $R(t)$ is $T$-periodic. They proved the existence of $T$-periodic solutions when $R$ is a positive $C^2$-function and $T<\pi$. When $R$ is a positive constant this equation is equivalent to \eqref{pin} and an isochronous center with minimal period $\pi$ appears. This suggests that the condition $T<\pi$ in the result in \cite{ACW} should be essential due to the appearance of resonance for $T=\pi$. We have analysed an example that somehow shows that this is the case, although our function $R(t)$ is not smooth. We thank Prof. Pedro J.~Torres for bringing to our attention the result in \cite{ACW}.

The rest of the paper is organized in five sections. The main results of the paper are stated in Section~\ref{sec:results}. Those readers who are only interested on the results for oscillators without singularities can go directly to Section~\ref{sec:resonance} to find the proofs. Section~\ref{sec:sturm} is concerned with some preliminary remarks on Sturm theory. In Section~\ref{sec:4} we present a detailed analysis of the asymptotic behaviour of the action-angle coordinates associated to~\eqref{aut}. These results will be employed in Section~\ref{sec:pinney} to prove the main theorems concerning oscillators with singularities. The paper is finished with an Appendix inspired by~\cite{BFS}.

\section{Statement of the main results}\label{sec:results}
Let us start with a potential $V\in C^2(\R)$ defined on the whole real line and satisfying
\[
V(0)=0,\ xV'(x)>0 \text{ if }x\neq 0.
\]
In addition we assume that all the solutions of equation \eqref{aut} are $2\pi$-periodic. In particular this implies that $V''(0)=N^2$, $N=1,2,\dots$ where $N$ depends upon the minimal period. See \cite{Urabe62} for more details and the construction of concrete examples.

There will be no particular restrictions on the forcing term $p(t)$ and we will just assume that it is $2\pi$-periodic and locally integrable. This will be indicated by $p\in L^1(\T)$.

\begin{thmx}\label{main}
Assume that $V$ satisfies the previous conditions and $V''$ is bounded over the whole real line. In addition the condition \eqref{rc} holds for some $p\in L^1(\T)$. Then the equation \eqref{per} is resonant for small $\epsilon\neq 0$.
\end{thmx}

The proof of this result will be presented in Section~\ref{sec:main}.

We observe that the potential associated to the asymmetric oscillator is not $C^2$. The reader is invited to modify the proof of the Theorem so that it includes this example.

In the introduction we defined resonant equation as an equation where all solutions are unbounded. This means that each solution $x(t)$ satisfies
\[
\abs{x(t_n)}+\abs{\dot x(t_n)}\rightarrow+\infty
\] 
for some sequence $\{t_n\}$. We will prove that a slightly stronger notion of resonance also holds, every solution $x(t)$ of \eqref{per} satisfies
\[
\lim_{\abs{t}\rightarrow+\infty}\left(\abs{x(t)}+\abs{\dot x(t)}\right)=+\infty.
\]

\begin{rem}\label{rem:def-global}
We point out that all the solutions of~\eqref{per} are globally defined. Indeed, if $x(t)$ is solution of~\eqref{per} and $E(t)=\frac{1}{2}\dot x(t)^2+V(x(t))$ then $E(t)$ satisfies the differential inequality 
\[
\abs{\dot E(t)}\leq \epsilon\abs{p(t)}\abs{\dot x(t)}\leq \epsilon\abs{p(t)}\sqrt{2E(t)}.\]
Thus, 
\begin{equation}\label{EEst}
\abs{\sqrt{E(t)}-\sqrt{E(0)}}\leq \frac{\epsilon}{\sqrt{2}}\abs{\int_0^t \abs{p(s)}ds}
\end{equation}
and so the energy $E(t)$ cannot blow up in finite time.
\end{rem}

Let us now assume that the function $p\in L^1(\T)$ is such that the condition \eqref{rc} does not hold. This means that $\Phi_p$ has a zero at some point $(\theta_*,r_*)$ of the cylinder $\mathcal C$ or that it vanishes at infinity. That is, either $\Phi_p(\theta_*,r_*)=0$ or $\Phi_p(\theta_n,r_n)\rightarrow 0$ for some sequence $(\theta_n,r_n)$ with $r_n\rightarrow+\infty$. The next result shows that resonance cannot occur if $\Phi_p$ has a non-degenerate zero.

\begin{prop}\label{reciprocal}
Assume that $V$ is in the conditions of Theorem~\ref{main} and that $\Phi_p$ has a non-degenerate zero $(\theta_*,r_*)$ with $r_*>0$. Then the equation \eqref{per} has a $2\pi$-periodic solution for small $\epsilon$.
\end{prop}

In the previous statement non-degeneracy is understood in the topological sense. That is, there exists a small open neighbourhood $\mathcal U$ of $(\theta_*,r_*)$ such that $\Phi_p(\theta,r)\neq 0$ for each point $(\theta,r)\in\bar{\mathcal U}\setminus\{(\theta_*,r_*)\}$ and the Brouwer degree $\deg(\Phi_p,\mathcal U,0)$ does not vanish. We refer to Section~\ref{sec:reciprocal} for the proof and more comments on Brouwer degree.

The two previous results show that the function $\Phi_p$ is crucial for the understanding of the resonance problem. In concrete examples it can be rather difficult to analyse the properties of $\Phi_p$. In the next result we show that, at least at a theoretical level, the condition \eqref{rc} is applicable to any isochronous potential and it is persistent under small perturbations of the forcing term $p(t)$. To be precise it is convenient to employ the language of Functional Analysis. Now we interpret $L^1(\T)$ as a Banach space with norm
\[
\|p\|_{L^1(\T)}=\int_0^{2\pi}\abs{p(t)}dt.
\]

\begin{prop}\label{functional}
Assume that $V$ is in the conditions of Theorem~\ref{main} and define $\mathcal R_V$ as the set of forcing terms satisfying \eqref{rc},
\[
\mathcal R_V=\{p\in L^1(\T) : \inf_{\mathcal C}\abs{\Phi_p(\theta,r)}>0\}.
\]
Then $\mathcal R_V$ is open and non-empty.
\end{prop}

This result will be proved in Section~\ref{sec:functional}. By now we observe that if $\{p_n\}$ is a sequence in $L^1(\T)$ converging in a weak topology to the Dirac measure $(p_n\wkconv*{}{}\delta)$ then $p_n\in\mathcal R_V$ for large $n$. This is consistent with the example in \cite{Ortega}. The convergence $p_n\wkconv*{}{}\delta$ means
\begin{enumerate}[$(a)$]
\item $\sup_n \|p_n\|_{L^1(\T)}<\infty,$ and
\item $\int_0^{2\pi} p_n(t)\phi(t)dt \rightarrow \phi(0)$ for each function $\phi(t)$ that is continuous and $2\pi$-periodic.
\end{enumerate}

In the previous results we have worked with oscillators defined on the whole real line but there are also oscillators producing an isochronous center and having a singularity. A well-known example is
\begin{equation}\label{pin}
\ddot x + \frac{1}{4}\left( x+1-\frac{1}{(x+1)^3}\right)=0,
\end{equation}
defined for all $x\in(-1,+\infty)$. This equation can be solved explicitly (see \cite{Pinney}). In particular
\[
\varphi(t,r)=-1+\sqrt{\lambda^2\cos^2(t/2)+\frac{1}{\lambda^2}\sin^2(t/2)}
\]
with $\lambda\!:=1+r\geq 1$. Let $\psi(t,r)$ be the solution of the associated variational equation
\begin{equation}\label{eq:var}
\ddot \psi + \frac{1}{4}\left(1+\frac{3}{(\varphi(t,r)+1)^4}\right)\psi =0
\end{equation}
with initial conditions $\psi(0,r)=1$, $\dot\psi(0,r)=i.$ A computation shows that
\begin{equation}\label{eq:expr-var}
\psi(t,r)=\frac{\partial \varphi}{\partial r}(t,r)-\frac{1}{V'(r)}\dot\varphi(t,r)i=\frac{\cos^2\bigl(\frac{t}{2}\bigr)-\frac{1}{\lambda^4}\sin^2\bigl(\frac{t}{2}\bigr)+2i\sin\bigl(\frac{t}{2}\bigr)\cos\bigl(\frac{t}{2}\bigr)}{\sqrt{\cos^2\bigl(\frac{t}{2}\bigr)+\frac{1}{\lambda^4}\sin^2\bigl(\frac{t}{2}\bigr)}}.
\end{equation}

The perturbed equation
\begin{equation}\label{bhora}
\ddot x + \frac{1}{4}\left(x+1-\frac{1}{(x+1)^3}\right)=\epsilon\sin t
\end{equation}
was considered in \cite{BFS}. The authors proved that all solutions are unbounded if $\epsilon\neq 0$ is small enough. Next we present an analogous result valid for general periodic perturbations of \eqref{pin}.

\begin{thmx}\label{thm:resonancia_pinney}
Let $p\in L^1(\T)$ be a function satisfying \eqref{rc} with $\psi(t,r)$ defined by \eqref{eq:expr-var}. Then all the solutions of equation
\begin{equation}\label{eq:per-pi}
\ddot x + \frac{1}{4}\left(x+1 - \frac{1}{(x+1)^3}\right) = \epsilon p(t)
\end{equation}
are unbounded for sufficiently small $\epsilon\neq 0$.
\end{thmx}

This resonance result deals with the specific equation \eqref{pin} but the method of proof can be extended to a larger class of potentials $V$. In many cases the existence of a singularity of $V$ at $x=-1$ determines the behaviour of $V(x)$ as $x\rightarrow+\infty$. This interesting observation was made in \cite{BFS} and somehow shows that the equation \eqref{pin} can be seen as a paradigm for centers with singularity. For more details we refer to the Appendix.

To show how to deduce the result for \eqref{bhora} in \cite{BFS} from Theorem~\ref{thm:resonancia_pinney} we present a corollary for the class of equations \eqref{eq:per-pi} with
\[
p(t)=a_0+a_1\cos t+b_1\sin t.
\]

\begin{cor}\label{corolario}
Assuming that $p(t)$ is a trigonometric function as above, the equation~\eqref{eq:per-pi} is resonant if 
\begin{equation}\label{New}
a_1^2+b_1^2>9a_0^2.
\end{equation}
\end{cor}

Note that this result extends the resonance result in~\cite{BFS} for the case $a_0=0$. 

We finish this Section with a result on the equation~\eqref{acw}.

\begin{thmx}\label{thm:multiplicativo}
Consider the $\pi$-periodic function
\[
R(t)=\begin{cases}
 1 & \text{ if } t\in\left[0,\frac{\pi}{2}\right),\\
 c & \text{ if } t\in\left[\frac{\pi}{2},\pi\right),
\end{cases}
\]
with $c>0$. Then all solutions of \eqref{acw} are unbounded if $c\neq 1$.
\end{thmx}

Since $R(t)$ is discontinuous the equation~\eqref{acw} is understood in the Carath\'eodory sense. It would be interesting to construct similar examples for smooth functions $R(t)$.

\section{Variations on Sturm Theory}\label{sec:sturm}
In the forthcoming sections we study the variations of the solution of the autonomous system~\eqref{aut} with respect to the action-angle variables. Some qualitative properties of these solutions will be understood thanks to a variant of Sturm theory. The classical results in this theory deal with the zeros of solutions but it is known that in some cases these zeros can be replaced by the zeros of the derivatives (critical points). See~\cite{Kamke}.

In this section we shall consider the equation
\begin{equation}\label{hill}
\ddot y + a(t)y=0
\end{equation}
with $a\in L^1(J)$ defined in an open interval $J$ and such that 
\begin{equation}\label{con}
\int_{J_*} a(t)dt>0 \ \text{ for each sub-interval }J_*\subset J.
\end{equation}
We know from Lebesgue differentiation theorem that
\[
\lim_{h\rightarrow 0} \frac{1}{h}\int_{t}^{t+h} a(s)ds = a(t)
\]
for almost every $t\in J$. In consequence $a(t)\geq 0$ almost everywhere.

Our basic tool will be the argument $\theta(t)$ of any non-trivial solution $y(t)$ of \eqref{hill}. Using the polar change of variables $y+i\dot y=re^{i\theta}$ we have that the argument satisfies
\begin{equation}\label{arg}
\dot\theta = -(\sin^2\theta + a(t)\cos^2\theta).
\end{equation}

Critical points of $y(t)$ correspond to $\theta\in\pi\Z$; that is, $\dot y(t)=0$ is equivalent to $\theta(t)=n\pi$ for some integer $n$. In the next result we prove that the trajectory $(y(t),\dot y(t))$ rotates around the origin in a clock-wise sense. This fact has useful consequences.

\begin{lema}\label{lema:argumento_monotono}
Let $y(t)$ be a non-trivial solution of $\eqref{hill}$. Then the argument $\theta(t)$ is strictly decreasing. In consequence,
\begin{enumerate}[$(a)$]
\item all critical points of $y(t)$ are isolated,
\item a local maximum or minimum is reached at each critical point of $y(t)$,
\item $y(t)$ reaches a local maximum (respectively, minimum) at $t^*\in J$ if and only if $\dot y(t^*)=0$, $y(t^*)>0$ (respectively, $y(t^*)<0$).
\end{enumerate}
\end{lema}

\begin{demo}
Given $t_1,t_2\in J$ with $t_1<t_2$ we integrate the equation~\eqref{arg} over the interval $J_*=(t_1,t_2)$ and observe that
\[
\theta(t_2)-\theta(t_1)=-\int_{t_1}^{t_2}(\sin^2\theta(t)+a(t)\cos^2\theta(t))dt\leq 
-\int_{t_1}^{t_2}\min\{1,a(t)\} dt<0.
\]
At this point we have employed the condition~\eqref{con}. Once we know that $(y(t),\dot y(t))$ rotates in a clock-wise sense we deduce that if $t^*\in J$ is such that $\theta(t^*)=n\pi$ with $n$ even, then there exists $\delta>0$ such that
\[
y(t)>0, \ (t-t^*)\dot y(t)<0\ \text{ if }0<\abs{t-t^*}<\delta.
\]
These inequalities are reversed when $n$ is odd. Now it is easy to prove all the properties $(a)-(c)$.
\end{demo}

Next we state a result on separation of critical points.

\begin{lema}\label{lema:maximos_alternados}
Let $\phi_1$ and $\phi_2$ be a fundamental pair of solutions of \eqref{hill}. 
Assume that $\phi_1$ reaches two local maxima (respectively, minima) at $t_1,t_2\in J$ with $t_1<t_2$. Then there exists $t_0\in(t_1,t_2)$ such that $\phi_2$ reaches a local maximum (respectively, minimum) at $t_0$.
\end{lema}

\begin{demo}
We assume that a local maximum is reached at $t_1$ and $t_2$. Otherwise we could replace $\phi_1$ by $-\phi_1$. Let $\theta_i(t)$ be the argument function corresponding to $\phi_i(t)$. Both arguments are solutions of the same first order equation~\eqref{arg} and so they cannot cross. Without loss of generality we can assume that $\theta_1(t_1)=0$ and $\theta_2(t_1)\in(-2\pi,0)$. Then $\theta_1(t)>\theta_2(t)$ whenever $t>t_1$. Since $\theta_1(t)$ is decreasing we know that $\theta_1(t_2)=-2n\pi$ for some integer $n\geq 1$. This implies that $\theta_2(t_2)<-2n\pi$ and so there exists some $t_0\in(t_1,t_2)$ such that $\theta_2(t_0)=-2\pi$. Consequently, a local maximum of $\phi_2$ is reached at $t_0$.
\end{demo}

We finish this section with a result on symmetric equations.

\begin{lema}\label{lema:2_criticos}
Let $\phi_1$ and $\phi_2$ be a fundamental pair of solutions of \eqref{hill}. Assume that $a(t)$ is defined in an open symmetric interval $J=(-\tau,\tau)$ such that $a(t)=a(-t)$ for all $t\in J$. Moreover, assume that $\dot\phi_1(0)=\phi_2(0)=0$ and $\phi_2$ has no critical points in $J$. Then $\phi_1$ has no critical points lying in $J\setminus\{0\}$.
\end{lema}

\begin{demo}
First we observe that the symmetry of the system with respect to $t=0$ induces the symmetries on the solutions $\phi_1(-t)=\phi_1(t)$ and $\phi_2(-t)=-\phi_2(t)$. The function $\phi_1(t)$ has a critical point at $t=0$, say that it is a maximum. Since $\phi_2$ has no critical point on $J$ then, by Lemma~\ref{lema:maximos_alternados}, $\phi_1$ has at most another critical point $t_*\neq 0$ (a minimum). The function $\phi_1$ is even and therefore critical points must be symmetric with respect to $t=0$. This excludes the possible existence of $t_*$ and so $t=0$ is the only critical point of $\phi_1$.
\end{demo}

\section{Potential centers and action-angle variables}\label{sec:4}

In this section we work with the equation
\begin{equation}\label{eq:system}
\ddot{x}+V'(x)=0,\ x\in\mathcal J
\end{equation}
where $\mathcal J=(a,b)$ is some interval with $-\infty\leq a<0<b\leq +\infty$ and the following hypothesis on the potential is considered
\begin{enumerate}
\item[$(H_0)$] $V\in\mathcal C^2(\mathcal J)$, $V(0)=V'(0)=0$, $V'(x)x>0$ if $x\neq 0$, $\lim_{x\rightarrow\partial\mathcal J}V(x)=+\infty$ and $V''(0)\neq 0$.
\end{enumerate}
In these conditions the equation~\eqref{eq:system} produces a center at the origin with $\mathcal J$ as a projection of the period annulus. Given a point $(x,\dot x)$ in the punctured strip $(\mathcal J\times \R)\setminus\{(0,0)\}$, there exists a unique closed orbit $\gamma$ passing through it. The action $I=I(x,\dot x)$ is the area of the region enclosed by $\gamma$. The angle $\theta=\theta(x,\dot x)$ can be defined mechanically as the quotient
\[
\theta = \frac{2\pi\tau}{T}
\]
where $\tau$ is the time employed to travel through $\gamma$ from the horizontal axis $\{\dot x=0\}$ to the point and $T$ is the minimal period of the orbit. Let $\gamma_r$ be the orbit passing through $(r,0)$ with $0<r<b$. The region enclosed by $\gamma_r$ will be denoted by $A_r$ and it increases with $r$. Moreover $\bigcup_{0<r<b}A_r=\mathcal J\times\R$ and so $I=\text{meas}(A_r)\rightarrow+\infty$ as $r\rightarrow b$. From these facts and the well-known theory of integrable systems we deduce that the map
\[
(x,\dot x)\in(\mathcal J\times\R)\setminus\{(0,0)\}\mapsto (I,\bar\theta)\in(0,\infty)\times(\R/2\pi\Z)
\]
defines a symplectic diffeomorphism transforming the equation~\eqref{eq:system} into
\[
\dot I = 0, \ \dot\theta = \omega(I)
\]
where $\omega(I)=\frac{2\pi}{T(I)}$ and $T(I)$ is the minimal period interpreted as a smooth function of the action.

Let $x(t,I)$ be the solution of~\eqref{eq:system} having action $I>0$ and satisfying the conditions $x(0,I)>0$, $\dot x(0,I)=0$. This section is concerned with the study of this function and the derivatives $\dot x=\frac{\partial x}{\partial t}$ and $\frac{\partial x}{\partial I}.$ Letting $r=x(0,I)$ we recall the classical formula
\begin{equation}\label{landau}
\Omega(I)=V(r)
\end{equation}
where $\Omega$ is the primitive of $\omega$ with $\Omega(0)=0$. This identity allows to connect the function $x(t,I)$ with the function $\varphi(t,r)$ appearing in the resonance condition~\eqref{rc}. Note that \eqref{landau} also shows that $\Omega(I)$ coincides with the energy along the corresponding orbit.

\begin{defi}\label{defi:P-N}
For each $I>0$ we define by 
\[
\mathcal P(I)\!:=\{t\in[0,T(I)]: x(t,I)\geq 0\}
\] 
the set of times when $x(t,I)$ is non-negative. We also define by $\mathcal N(I)\!:=[0,T(I)]\setminus\mathcal P(I)$ the set of times when $x(t,I)$ is negative. 
\end{defi}

Notice that, due to the symmetry $x(-t,I)=x(t,I)$, $\mathcal N(I)$ is an interval centered at $t=\frac{T(I)}{2}$ and $\mathcal P(I)$ is the union of two intervals.

\begin{lema}\label{lema:x_partial_I}
The function $x(t,I)$ belongs to $C^1(\R\times(0,\infty))$ with
\begin{equation}\label{formula_xI}
\frac{1}{\omega(I)}\frac{\partial x}{\partial I}(t,I)=-\frac{\ddot{x}(t,I)}{\dot{x}(t,I)^2+\ddot{x}(t,I)^2}+\dot{x}(t,I)\int_0^t\frac{\left(1-V''\bigl(x(s,I)\bigr)\right)(\dot{x}(s,I)^2-\ddot{x}(s,I)^2)}{\bigl(\dot{x}(s,I)^2+\ddot{x}(s,I)^2\bigr)^2}ds.
\end{equation}
Moreover, $\frac{\partial x}{\partial I}\bigl(\frac{T(I)}{2},I\bigr)<0$.
\end{lema}

\begin{demo}
The smoothness of $x(t,I)$ follows from the connection with $\varphi(t,r)$ via the identity~\eqref{landau}. Indeed the functions $\dot x(t,I)$ and $\frac{\partial x}{\partial I}(t,I)$ are solutions of the variational equation of~\eqref{eq:system} given by
\begin{equation}\label{eq:variational}
\ddot{\xi}+V''\bigl(x(t,I)\bigr)\xi=0.
\end{equation}
From~\eqref{landau} we deduce that $\frac{\partial x}{\partial I}(t,I)$ satisfies the initial conditions
\[
\frac{\partial x}{\partial I}(0,I)=\frac{\omega(I)}{V'(r)}, \ \frac{\partial \dot x}{\partial I}(0,I)=0.
\]
The Rofe-Beketov formula (see \cite{BES}, page $24$) implies that the function $\xi(t)$ given by the right hand side of~\eqref{formula_xI} is a solution of~\eqref{eq:variational} and a direct computation shows that
\[
\xi(0)=-\frac{1}{\ddot x(0,I)}=\frac{1}{V'(r)},\ \dot\xi(0)=0.
\]
The uniqueness of solution of the initial value problem implies that the identity~\eqref{formula_xI} holds.
The inequality $\frac{\partial x}{\partial I}\bigl(\frac{T(I)}{2},I\bigr)<0$ follows from the expression of $\frac{\partial x}{\partial I}(t,I)$ with $t=\frac{T(I)}{2}$ using that $\dot{x}\bigl(\frac{T(I)}{2},I\bigr)=0$ and $\ddot{x}\bigl(\frac{T(I)}{2},I\bigr)=-V'\bigl(x\bigl(\frac{T(I)}{2},I\bigr)\bigr)>0$.
\end{demo}

\begin{lema}\label{lema:bounds}
In the previous notations, $\abs{\dot{x}(t,I)}\leq \sqrt{2\Omega(I)}$ for all $t\in\R$ and $I>0$.
\end{lema}

\begin{demo}
By energy conservation we have $\frac{1}{2}\dot{x}(t,I)^2=\Omega(I)-V\bigl(x(t,I)\bigr)$. Then, taking into account that $V(x)\geq 0$ the result holds.
\end{demo}

The following sections are concerned with the behaviour of the periodic solutions $x(t,I)$ of system \eqref{eq:system} and its partial derivatives with respect to the action-angle variables. Section~\ref{sec:asintota} is devoted to potential centers that are not globally defined by an asymptote of the potential function. In Section~\ref{sec:lipschitz} we consider the case when the second derivative of the potential is bounded. Section~\ref{sec:ambos} is concerned with a combination of the previous situations, an asymptote on one side and a bounded second derivative of $V$ on the other side. Finally, in Section~\ref{sec:comportamiento} we continue with this situation but in addition we assume that the center is isochronous. We describe the asymptotic behaviour of $x(t,I)$, $\dot x(t,I)$ and $\frac{\partial x}{\partial I}(t,I)$ as $I\rightarrow+\infty$. The limit functions are solutions of a Bouncing Problem that will be analyzed in detail.

\subsection{Potential center with an asymptote}\label{sec:asintota}

In this section we shall consider that the potential function $V$ in \eqref{eq:system} presents an asymptote at $x=a$. That is, we shall assume that $V$ satisfies $(H_0)$ with $\mathcal J=(a,+\infty)$, $-\infty<a<0$. In addition we shall assume that $V$ is convex near the asymptote, meaning that 
\begin{enumerate}
\item[$(H_1)$] $V''(x)>0$ for all $x\in(a,a+\zeta)$.
\end{enumerate}
where $\zeta$ is some number satisfying $0<\zeta<\abs{a}$.

The first result is concerned with the limit as $I\rightarrow+\infty$ of the negative semi-period of $x(t,I)$. It is based on a similar result in~\cite{ManRojVil2015}. In order to state it properly, denote by $T_-(I)$ the negative semi-period for each $I\in(0,+\infty)$. That is, the length of the interval $\mathcal N(I)$. Note that
\[
\mathcal N(I)=\left(\frac{T(I)}{2}-\frac{T_-(I)}{2},\frac{T(I)}{2}+\frac{T_-(I)}{2}\right).
\]

\begin{lema}\label{lema:semiperiodo-cero}
Assume that $V$ satisfies $(H_0)$ and $(H_1)$. Then $T_-(I)\rightarrow 0$ as $I\rightarrow+\infty$.
\end{lema}

In the proof of this result we will employ the following elementary property.

\begin{lema}\label{lema:convex}
Let $f:[M,+\infty)\rightarrow\R$ be a convex and decreasing function. In addition assume that $f$ has derivative $f'(x)$ everywhere and $\lim_{x\rightarrow+\infty}f(x)$ exists and is finite. Then 
\[
\lim_{x\rightarrow+\infty}xf'(x)=0.
\]
\end{lema}

\begin{demo}
Assume $M>0$. By the convexity property we have that for all $x,y\in[M,+\infty)$, $f(y)\geq f(x)+f'(x)(y-x)$. Taking $x>2M$ and $y=x/2$, the previous inequality yields to $xf'(x)\geq 2(f(x)-f(x/2))$ and, using that $f$ is decreasing, $\abs{xf'(x)}\leq 2\abs{f(x)-f(x/2)}$ for all $x>2M$. Then, since $f(x)$ has a finite limit as $x\rightarrow+\infty$ the result holds.
\end{demo}

\begin{prooftext}{Proof of Lemma~\ref{lema:semiperiodo-cero}}
We denote by $\gamma_-(I)$ the portion of the orbit described by $x(t,I)$ on the half negative plane and by $V_-$ the restriction of the potential $V$ on the negative semi-axis. Then, since $\Omega(I)$ is the energy of the orbit, the negative semi-period is written as
\[
T_-(I)=\int_{\gamma_-(I)}\frac{dx}{y}=\sqrt{2}\int_{r_-(I)}^{0}\frac{dx}{\sqrt{\Omega(I)-V_-(x)}},
\]
where $r_{-}(I)$ denotes the intersection of $\gamma_-(I)$ with the $x$-axis.
The change of variable $V_-(x)^{\frac{1}{2}}=-\Omega(I)^{\frac{1}{2}}\sin\theta$ transforms the previous integral into
\[
T_-(I)=-2\sqrt{2}\int_{-\frac{\pi}{2}}^0\left(\frac{(V_-)^{\frac{1}{2}}}{V_-'}\right)\bigl(V_-^{-1}(\Omega(I)\sin^2\theta)\bigr)d\theta.
\]
From hypothesis $(H_0)$ we have that  
\begin{equation}\label{eq:limite_cero}
\lim_{x\rightarrow 0^-}\frac{(V_-)^{\frac{1}{2}}(x)}{V_-'(x)}= - \frac{1}{\sqrt{2V''(0)}}.
\end{equation}
On the other hand, by hypothesis $(H_1)$, $V_-^{-1}$ is a convex and decreasing function defined in $(V_-(a+\zeta),+\infty)$ with limit $a<0$ at infinity so, by Lemma~\ref{lema:convex} and setting $y=V_-(x)$, we have that
\begin{equation}\label{eq:limite_a}
\lim_{x\rightarrow a}\frac{(V_-)^{\frac{1}{2}}(x)}{V_-'(x)}=\lim_{y\rightarrow+\infty}\frac{y(V_-^{-1})'(y)}{y^{\frac{1}{2}}}=0.
\end{equation}
Therefore, on account of \eqref{eq:limite_cero} and \eqref{eq:limite_a}, there exists $M>0$ such that $\abs{\frac{(V_-)^{\frac{1}{2}}}{V_-'}\circ V_-^{-1}(y)}\leq M$ for all $y\in(0,+\infty)$. Since $\Omega(I)\rightarrow+\infty$ as $I\rightarrow+\infty$ then $V_-^{-1}(\Omega(I)\sin^2\theta)\rightarrow a$ for every $\theta\in[-\frac{\pi}{2},0)$.
Thus, by Lebesgue's dominated convergence theorem and on account of \eqref{eq:limite_a}, 
$T_-(I)\rightarrow 0$ as $I\rightarrow+\infty$.
\end{prooftext}

The previous result states that the presence of an asymptote of the potential at the negative $x$-axis implies that the periodic solution $x(t,I)$ tends to stay on the positive half-plane as the action tends to infinity. 

A significant role will be played by the derivative of the periodic solutions with respect to the action on the interval $\mathcal N(I)$. The next result will help to control its behaviour. Before stating it we introduce some more notation.

\begin{defi}\label{defi:A}
Let $I_{\zeta}$ be defined by the equation $\Omega(I_{\zeta})=V(a+\zeta)$. This is the value of the action corresponding to the energy at the point $(a+\zeta,0)$. For each $I>0$ we define
\[
\mathcal A(I)\!:=(\tau^-(I),\tau^+(I))=\{t\in[0,T(I)]: x(t,I)<a+\zeta\}
\] 
the interval of time when $x(t,I)$ belongs to the interval $(a,a+\zeta)$.
\end{defi}
Due to the symmetry $x\bigl(t+\frac{T(I)}{2},I\bigr)=x\bigl(-t+\frac{T(I)}{2},I\bigr)$, when $I\geq I_{\zeta}$ then $\mathcal A(I)$ is an interval centered at $t=\frac{T(I)}{2}$ and $\mathcal A(I)=\{\frac{T(I)}{2}\}$ if and only if $I=I_{\zeta}$. Otherwise $\mathcal A(I)$ is empty.

\begin{lema}\label{lema:cota_xI_pi}
Let $x(t,I)$ be a solution of system \eqref{eq:system} satisfying $(H_0)$ and $(H_1)$. Then the following holds:
\begin{enumerate}[$(a)$]
\item $\dfrac{\partial x}{\partial I}(T(I)/2,I)\leq \dfrac{\partial x}{\partial I}(t,I) \leq \dfrac{\partial x}{\partial I}(\tau^+(I),I) \text{ for all }t\in \mathcal A(I).$
\item $0<-\dfrac{\partial x}{\partial I}(T(I)/2,I)\leq \abs{a}\dfrac{\omega(I)}{I-I_{\zeta}}$ for all $I>I_{\zeta}$.
\end{enumerate} 
\end{lema}

\begin{demo}
The functions $\dot{x}(t,I)$ and $\frac{\partial x}{\partial I}(t,I)$ are a fundamental pair of solutions of the variational equation $\ddot \xi + V''(x(t,I))\xi =0$. By definition of $\mathcal A(I)$, we have that $V''(x(t,I))>0$ for all $t\in\mathcal A(I)$ and $\dot{x}(t,I)$ has no critical points for $t\in\mathcal A(I)$. Therefore the variational equation satisfies the condition~\eqref{con} on $\mathcal A(I)$ and, by Lemma~\ref{lema:2_criticos}, the only critical point of 
$\frac{\partial x}{\partial I}(t,I)$ in $\mathcal A(I)$ is $t=\frac{T(I)}{2}$. 
Moreover, by Lemma~\ref{lema:x_partial_I} we have $\frac{\partial x}{\partial I}\bigl(\frac{T(I)}{2},I\bigr)<0$ so $t=\frac{T(I)}{2}$ is a minimum as a consequence of Lemma~\ref{lema:argumento_monotono} and so $(a)$ is proved.

To prove $(b)$ let us denote $r_{-}(I)\!:=x\bigl(\frac{T(I)}{2},I\bigr)$. We have then $r_{-}'(I)=\frac{\partial x}{\partial I}\bigl(\frac{T(I)}{2},I\bigr)$ due to $\dot x\bigl(\frac{T(I)}{2},I\bigr)=0$. On the other hand, by conservation of energy $V(r_{-}(I))=\Omega(I)$ and so $\frac{\partial x}{\partial I}\bigl(\frac{T(I)}{2},I\bigr)=\frac{\omega(I)}{V'(r_{-}(I))}$. Since $I>I_{\zeta}$ we have $a<r_{-}(I)<a+\zeta$. Consequently, there exists $\eta=\eta(I)\in(r_-(I),a+\zeta)$ such that 
\[
I-I_{\zeta}=V(r_{-}(I))-V(a+\zeta)=V'(\eta)(r_{-}(I)-a-\zeta)
\]
for all $I>I_{\zeta}$. On account of hypothesis $(H_1)$, $\abs{V'(r_{-}(I))}\geq \abs{V'(\eta)}$ for all $I>I_{\zeta}$. Hence,
\[
\abs{\frac{\partial x}{\partial I}(T(I)/2,I)}\leq \frac{\omega(I)}{\abs{V'(\eta)}}\leq \omega(I)\frac{a+\zeta-r_{-}(I)}{I-I_{\zeta}}\leq  \omega(I)\frac{\abs{a}}{I-I_{\zeta}}
\] 
for all $I>I_{\zeta}$ as we desired.
\end{demo}

\begin{rem}
In the above results the potential function $V(x)$ has an asymptote at $x=a$, $a<0$. Analogous results can be obtained when $a=-\infty$ and the asymptote is located at $b\in(0,+\infty)$. Indeed it is enough to use the change of variables $x\mapsto -x$.
\end{rem}

\subsection{Potential center with Lipschitz behaviour at infinity}\label{sec:lipschitz}
In this section we consider a potential function defined on $\mathcal J=(a,+\infty)$ with $-\infty\leq a <0$ and satisfying $(H_0)$. In addition we assume that
\begin{enumerate}
\item[$(H_2)$] $V''$ is bounded in $J$
\end{enumerate}
where $J$ is a fixed interval of type $J=[\alpha,+\infty)$ with $a<\alpha\leq 0$.
We denote $\|V''\|_{\infty}\!:=\sup\{\abs{V''(x)}:x\in J\}<\infty$.
\begin{defi}\label{defi:B}
Let $x(t,I)$ be a solution of system \eqref{eq:system} satisfying hypothesis $(H_0)$ and $(H_2)$. For each $I>0$, we define  
\[
\mathcal B(I)\!:=\{t\in[0,T(I)]:x(t,I)\in J\}.
\]
\end{defi}
We observe that either $\mathcal B(I)$ is the whole interval $[0,T(I)]$ or it is the union of two intervals containing respectively $t=0$ and $t=T(I)$.

\begin{lema}\label{lema:bounds2}
Let $x(t,I)$ be a solution of system \eqref{eq:system} satisfying $(H_0)$ and $(H_2)$. Then 
\[
\dot{x}(t,I)^2+\ddot{x}(t,I)^2\geq \frac{1}{4}e^{-2T(I)M}\ddot{x}(0,I)^2>0
\]
for all $t\in\mathcal B(I)$, where $M\!:=\max\{1,\|V''\|_{\infty}\}$. Note that $\abs{\ddot x(0,I)}=V'(r)$.
\end{lema}

\begin{demo}
The identity $x(T(I)-t,I)=x(t,I)$ implies that it is sufficient to prove the estimate on the component of $\mathcal B(I)$ containing $t=0$. Let us consider the variational equation $\ddot{\xi}+V''\bigl(x(t,I)\bigr)\xi=0$. Clearly $\dot{x}(t,I)$ is a solution of this equation. If we set $\eta\!:=\dot{\xi}$, $\mathcal X\!:=(\xi,\eta)^*$ and 
\[
A(t,I)\!:=\left(\begin{matrix} 0 & 1\\ -V''(x(t,I)) & 0 \\
\end{matrix}\right),
\] 
the variational equation is written as $\dot{\mathcal X}=A(t,I)\mathcal X$. Let $\Phi(t,I)$ be the fundamental matrix of the variational equation with $\Phi(0,I)=id$ for all $I>0$. The Gronwall's inequality states that 
\[
\|\Phi(t,I)\|\leq \exp\left(\int_{0}^t \|A(s,I)\|ds\right),
\]
where $\|\cdot\|$ denotes the matrix norm induced by the sup norm of $\R^2$. Given a matrix $A=(a_{ij})\in M_2(\R)$, $\|A\|=\max\{\abs{a_{11}}+\abs{a_{12}},\abs{a_{21}}+\abs{a_{22}}\}$. Since $t\in\mathcal B(I)$ we have $\|A(t,I)\|=\max\{1,\|V''\|_{\infty}\}=M$ and so $\|\Phi(t,I)\|\leq e^{T(I)M}$ for all $t\in\mathcal B(I)$. 
On the other hand, since the matrix $\Phi(t,I)$ is symplectic, $\Phi(t,I)^{-1}=-J\Phi(t,I)^*J$, where 
$J=\bigl(\begin{smallmatrix}
0 & 1\\ -1 & 0
\end{smallmatrix}\bigr)$.
Consequently, $\|\Phi(t,I)^{-1}\|\leq \|\Phi(t,I)^*\|\leq 2\|\Phi(t,I)\|$ and so
\[
\|\mathcal X(0)\|\leq \|\Phi(t,I)^{-1}\|\|\mathcal X(t)\|\leq 2e^{T(I)M}\|\mathcal X(t)\|\] 
for all $t$ lying in the component of $\mathcal B(I)$ containing $t=0$.
\end{demo}

\subsection{Potential center with an asymptote and Lipschitz behaviour at infinity}\label{sec:ambos}

In this section we consider potential functions $V\in C^2(\mathcal J)$, $\mathcal J=(a,+\infty)$ with $-\infty<a<0$ satisfying the hypotheses $(H_0)$, $(H_1)$ and $(H_2)$. Moreover,
\begin{enumerate}
\item[$(H_3)$] $\displaystyle\liminf_{x\rightarrow+\infty} \frac{V'(x)}{x}>0$ and $V''(0)>0$.
\end{enumerate}

In the above conditions it is not restrictive to assume that $a+\zeta=\alpha$. From now on, $\|V''\|_{\infty}\!:=\sup\{\abs{V''(x)}:x>a+\zeta\}$.

\begin{lema}\label{lema:amplitud}
Assume that $(H_i)$, $i=0,1,2,3,$ hold. Then there exists $c>0$ such that
\[
\frac{c}{2}r^2 \leq \Omega(I) \leq \frac{\|V''\|_{\infty}}{2}r^2, \ I>0.
\]
Here $x(t,I)$ is the solution of~\eqref{eq:system} and $r=x(0,I)$.
\end{lema}

\begin{demo}
The function $W(x)=\frac{V'(x)}{x}$, $x\in(0,+\infty)$ is positive and satisfies $W(x)\rightarrow V''(0)>0$ as $x\rightarrow 0^+$. Since $\liminf_{x\rightarrow+\infty} W(x)>0$ we conclude that $W(x)$ has a positive lower bound, say $c>0$. Once we know that
\begin{equation}\label{DV}
V'(x)\geq cx \text{ if }x>0
\end{equation}
we deduce that $V(x)\geq \frac{c}{2}x^2$ if $x>0$. Here we have used that $V(0)=0$. We can invoke the identity~\eqref{landau} to deduce the inequality $\Omega(I)=V(r)\geq \frac{c}{2}r^2$. On the other hand, by $(H_2)$ we have that $V''(x)\leq \|V''\|_{\infty}$ for all $x\geq 0$. Then, taking into account that $V(0)=V'(0)=0$,
\[
V(x)=\int_0^x\int_0^s V''(u)duds\leq \frac{\|V''\|_{\infty}}{2}x^2.
\]
The inequality $\Omega(I)\leq \frac{\|V''\|_{\infty}}{2}r^2$ follows.
\end{demo}

We are now in position to prove the main result of this section. This result gives upper-bounds of $x(t,I)$ and its partial derivatives in terms of the action. Since we are assuming $\alpha=a+\zeta$, the interval $[0,T(I)]$ splits as $\mathcal A(I)\cup\mathcal B(I)$.

\begin{prop}\label{prop:cotas}
Let $x(t,I)$ be a solution of system \eqref{eq:system} satisfying $(H_i)$, $i=0,1,2,3$. Then the following holds:
\begin{enumerate}[$(a)$]
\item $\abs{\dfrac{1}{\sqrt{\Omega(I)}}\ x(t,I)}\leq 
\begin{cases}
 \sqrt{\dfrac{2}{c}} & \text{ for all }t\in\mathcal P(I) \text{ and } I>0,\\
\dfrac{\abs{a}}{\sqrt{\Omega(I)} }& \text{ for all }t\in\mathcal N(I) \text{ and } I>0.
\end{cases}$
\item $\abs{\dfrac{1}{\sqrt{\Omega(I)}}\ \dot{x}(t,I)}\leq \sqrt{2}$ for all $t\in[0,T(I)]$ and $I>0$.
\item $
\abs{\sqrt{\Omega(I)}\ \dfrac{\partial x}{\partial I}(t,I)}\leq 
\begin{cases}
\Delta(I) & \text{ if }t\in\mathcal B(I),\\
\max\left\{\Delta(I),\dfrac{\abs{a}\omega(I)\sqrt{\Omega(I)}}{I-I_{\zeta}}\right\} & \text{ if }t\in \mathcal A(I),\ I>I_{\zeta}
\end{cases}
$
\end{enumerate}
where $\Delta(I)\!:=\frac{\sqrt{2M}}{c}e^{T(I) M} + \frac{4\sqrt{2}M^2T(I)}{c^2}e^{2T(I)M}$  and $M\!:=\max\{1,\|V''\|_{\infty}\}$.
\end{prop}

\begin{demo}
In order to prove $(a)$ let us note that by Lemma~\ref{lema:amplitud}, $\sqrt{\Omega(I)}\geq \sqrt{\frac{c}{2}}r$. Since $\abs{x(t,I)}\leq r$ for all $t\in\mathcal P(I)$ and $I>0$, we have 
\[
\abs{\frac{1}{\sqrt{\Omega(I)}}\ x(t,I)}\leq \sqrt{\frac{2}{c}}. 
\]
On the other hand, if $t\in\mathcal N(I)$ then $a<x(t,I)<0$ so the assertion in $(a)$ holds.

The proof of $(b)$ follows straightforward from Lemma~\ref{lema:bounds}. For proving $(c)$, let us first start by considering $t\in\mathcal B(I)$. Indeed it is not restrictive to assume that $t$ lies on the component containing $0$. By Lemma~\ref{lema:x_partial_I} we have
\[
\abs{\frac{\partial x}{\partial I}(t,I)}\leq 
\frac{\abs{\ddot{x}(t,I)}}{\dot{x}(t,I)^2+\ddot{x}(t,I)^2}+\abs{\dot{x}(t,I)}\int_0^t\frac{\abs{1-V''\bigl(x(s,I)\bigr)}\abs{\dot{x}(s,I)^2-\ddot{x}(s,I)^2}}{\bigl(\dot{x}(s,I)^2+\ddot{x}(s,I)^2\bigr)^2}ds.
\]
By hypothesis $(H_2)$ in $(a+\zeta,+\infty)$ we have that $\mathcal B(I)=[0,T(I)]\setminus\mathcal A(I)$ correspond to the set of times when $V''(x(t,I))$ is bounded. Therefore, on account of Lemma~\ref{lema:bounds} and Lemma~\ref{lema:bounds2} together with the general inequalities $a\leq\sqrt{a^2+b^2}$ and $a^2-b^2\leq a^2+b^2$ we have
\begin{align*}
\abs{\frac{\partial x}{\partial I}(t,I)}&\leq 
\frac{1}{(\dot x(t,I)^2+\ddot x(t,I)^2)^{\frac{1}{2}}}+\abs{\dot x(t,I)}\int_0^t\frac{\abs{1-V''(x(s,I))}}{\dot x(s,I)^2+\ddot x(s,I)^2}ds\\
&\leq
\frac{2e^{T(I)M}}{V'(r)}+\frac{4\sqrt{2\Omega(I)}e^{2T(I)M}}{V'(r)^2}\int_0^t\abs{1-V''(x(s,I))}ds.
\end{align*}
From \eqref{DV} we deduce that $V'(r)\geq cr\geq c\sqrt{\frac{2}{\|V''\|_{\infty}}}\Omega(I)^{\frac{1}{2}}$ leading to
\begin{equation}\label{eq:cota_en_P}
\abs{\sqrt{\Omega(I)}\frac{\partial x}{\partial I}(t,I)}\leq 
\frac{\sqrt{2M}}{c}e^{T(I) M} + \frac{4\sqrt{2}M^2T(I)}{c^2}e^{2T(I)M}
\end{equation}
for all $t\in\mathcal B(I)$. Here we also used that $1+\|V''\|_{\infty}\leq 2M$.

On the other hand, if $t\in\mathcal A(I)$, by $(a)$ in Lemma~\ref{lema:cota_xI_pi} we have
\[
\abs{\sqrt{\Omega(I)}\ \frac{\partial x}{\partial I}(t,I)}\leq \max\left\{\left|\sqrt{\Omega(I)}\ \frac{\partial x}{\partial I}(\pi,I)\right|,\left|\sqrt{\Omega(I)}\ \frac{\partial x}{\partial I}(\tau,I)\right|\right\},
\]
where $\tau=\tau_{\zeta}^+(I)$ is the right-hand endpoint of $\mathcal A(I)$. Notice that $\tau\in\mathcal B(I)$, so the second function in the previous maximum satisfies the inequality in \eqref{eq:cota_en_P} and so the result follows by $(b)$ in Lemma~\ref{lema:cota_xI_pi}.
\end{demo}

\subsection{Isochronous oscillators with an asymptote and a bouncing problem}\label{sec:comportamiento}

In this section we assume that $-\infty<a<0<b=+\infty$ and the potential $V$ satisfies $(H_i)$, $i=0,1,2,3$. In addition it will be assumed that all non-constant solutions of equation~\eqref{eq:system} have minimal period $2\pi$, that is $\Omega(I)=I$.
We will be interested in the asymptotic behaviour of the functions $x(t,I)$, $\dot x(t,I)$ and $\frac{\partial x}{\partial I}(t,I)$ as $I\rightarrow+\infty$. The following notion will be useful to describe this behaviour.

\begin{defi}\label{defi:BP}
Given $\gamma\in L^{\infty}(\T)$, we say that $\psi$ is a solution of the $\gamma-$\emph{Bouncing problem} $(\gamma-BP)$ if it satisfies the following properties:
\begin{enumerate}[$(i)$]
\item $\psi\in W^{1,\infty}(\T)\cap W^{2,\infty}(-\pi,\pi)$,
\item $\ddot\psi+\gamma(t)\psi=0$ in $(-\pi,\pi)$,
\item $\psi\geq 0$ and $\psi(t)=0$ if and only if $t=(2n+1)\pi$,
\item $\dot\psi(\pi+0)=-\dot\psi(\pi-0)$.
\end{enumerate}\vspace*{-.5cm}
\end{defi}

The space $L^{\infty}(\T)$ is composed by the measurable functions which are $2\pi$-periodic and essentially bounded. $W^{1,\infty}(\T)$ is the Sobolev space modelled on $L^{\infty}(\T)$. Indeed condition $(i)$ can be formulated in more classical terms, without reference to Sobolev spaces. More precisely, we can replace it by
\begin{enumerate}[$(i')$]
\item $\psi:\R\rightarrow\R$ is continuous and $2\pi$-periodic, the restriction $\left.\psi\right|_{[-\pi,\pi]}$ belongs to $C^1[-\pi,\pi]$ and the derivative $\psi'$ is Lipschitz-continuous on $[-\pi,\pi]$.
\end{enumerate}
The equation in $(ii)$ can be understood in the Carath\'eodory sense. Finally we observe that, by periodicity, $\psi$ also belongs to $C^1[\pi,3\pi]$ and so the condition $(iv)$ has a precise meaning. Typically $\psi'$ will present a jump discontinuity at the points $t=(2n+1)\pi$, $n\in\Z$.

\begin{thm}\label{thm:limite-x}
Let $(I_n)_{n\geq 0}$ be a sequence of actions such that $I_n\rightarrow+\infty$ as $n\rightarrow+\infty$. Then there exist a subsequence $(I_{n_k})_{k\geq 0}\subset(I_n)_{n\geq 0}$, $\gamma\in L^{\infty}(\T)$ and a solution $\psi$ of the $\gamma-BP$ such that
\begin{enumerate}[$(a)$]
\item $\dfrac{x(t,I_{n_k})}{\sqrt{I_{n_k}}}$ converges uniformly to $\psi$, and
\item $\dfrac{\dot x(t,I_{n_k})}{\sqrt{I_{n_k}}}$ converges uniformly to $\dot\psi$ in $[-\pi+\delta,\pi-\delta]$ for all $\delta>0$.
\end{enumerate}
Moreover the function $\gamma$ satisfies 
\begin{equation}\label{lab}
\liminf_{x\rightarrow+\infty}\frac{V'(x)}{x}\leq \gamma(t)\leq \limsup_{x\rightarrow+\infty}\frac{V'(x)}{x}
\end{equation}
almost everywhere.
\end{thm}

An analogous result holds for $\frac{\partial x}{\partial I}$.

\begin{thm}\label{thm:limite-xI}
In the same conditions of the previous theorem there exists a subsequence $(I_{m_k})_{k\geq 0}\subset(I_n)_{n\geq 0}$, $\gamma_1\in L^{\infty}(\T)$ and a solution $\psi_1$ of the $\gamma_1-BP$ such that
\[
\sqrt{I_{m_k}}\dfrac{\partial x}{\partial I}(t,I_{m_k}) \text{ converges uniformly to }\psi_1 \text{ in }[-\pi+\delta,\pi-\delta] \text{ for all }\delta>0.
\]
The function $\gamma_1$ satisfies the condition
\begin{equation}\label{lHop}
\liminf_{x\rightarrow+\infty}V''(x)\leq \gamma_1(t)\leq \limsup_{x\rightarrow+\infty} V''(x)
\end{equation}
almost everywhere.
\end{thm}

In view of the generalized L'H\^opital rule this condition is weaker than~\eqref{lab}. It must be noticed that the functions $\gamma$ and $\gamma_1$ do not necessarily coincide even when a common subsequence of $(I_n)$ can be extracted.

Before proving these results we need some preliminary lemmas about the weak$^*$ convergence in $L^{\infty}(J)$, where $J$ is a bounded interval. We first recall this notion of convergence. Given a sequence of functions $\phi_n\in L^{\infty}(J)$ and $\phi\in L^{\infty}(J)$, we say that $\phi_n$ converges to $\phi$ in the weak$^*$ sense, $\phi_n\wkconv*{*}{}\phi$ in $L^{\infty}(J)$, if
\[
\int_J \chi \phi_n\rightarrow\int_J \chi\phi \text{ for each }\chi\in L^1(J).
\]

\begin{lema}\label{lema:weak}
Let $(\phi_n)_{n\geq 0}$ be a sequence of functions, $\phi_n\in L^{\infty}(0,\pi)$, and $C>0$ such that for each $\delta>0$ there exists $N=N(\delta)$ satisfying $\|\phi_n\|_{L^{\infty}(0,\pi-\delta)}\leq C$ if $n\geq N$. Then there exist a subsequence $(\phi_{n_k})_{k\geq 0}\subset(\phi_n)_{n\geq 0}$ and $\hat\phi\in L^{\infty}(0,\pi)$ such that $\phi_{n_k}|_{(0,\pi-\delta)}\wkconv*{*}{}\hat\phi|_{(0,\pi-\delta)}$ in $L^{\infty}(0,\pi-\delta)$ for all $\delta>0$.
\end{lema}

\begin{demo}
Let $(\delta_k)_{k\geq 0}$ be a monotone sequence tending to zero and let us denote $X_k= L^{\infty}(0,\pi-\delta_k)$. After extending by zero any function of $X_k$ to $(0,\pi)$ we have
\[
X_0\subset X_1 \subset \cdots \subset X_k \subset \cdots \subset L^{\infty}(0,\pi).
\]
We will employ a diagonal procedure. The sequence $\|\phi_n\|_{X_0}$ is eventually bounded by $C>0$ by hypothesis. Thus, Banach-Alaoglu Theorem states that there exist $\hat \phi_0\in X_0$ and a subsequence $(\phi_{\sigma_0(n)})_{n\geq 0}\subset (\phi_n)_{n\geq 0}$ such that $\phi_{\sigma_0(n)}\wkconv*{*}{}\hat\phi_0$ in $X_0$. Moreover, $\|\hat\phi_0\|_{X_0}\leq C$. The sequence $\|\phi_n\|_{X_1}$ is also eventually bounded by $C>0$ by hypothesis. Thus, again by Banach-Alaoglu Theorem, there exist $\hat \phi_1\in X_1$ and a subsequence $(\phi_{\sigma_1(n)})_{n\geq 0}\subset (\phi_{\sigma_0(n)})_{n\geq 0}$ such that $\phi_{\sigma_1(n)}\wkconv*{*}{}\hat\phi_1$ in $X_1$. Moreover, $\|\hat\phi_1\|_{X_1}\leq C$ and $\hat\phi_1=\hat\phi_0$ almost everywhere in $(0,\pi-\delta_0)$. Here we are using that $(\phi_{\sigma_1(n)})$ is a subsequence of $(\phi_{\sigma_0(n)})$ and the uniqueness of the limit in the weak$^*$ sense. Inductively we have that for each $k\geq 0$ there exist $\hat\phi_k\in X_k$ and a subsequence $(\phi_{\sigma_k(n)})_{n\geq 0}\subset (\phi_{\sigma_{k-1}(n)})_{n\geq 0}\subset \cdots \subset (\phi_{\sigma_0(n)})_{n\geq 0}\subset (\phi_n)_{n\geq 0}$ such that $\phi_{\sigma_k(n)}\wkconv*{*}{}\hat\phi_k$ in $X_k$, $\|\hat\phi_k\|_{X_k}\leq C$ and $\hat\phi_k=\hat\phi_j$ almost everywhere in $(0,\pi-\delta_j)$ for all $0\leq j<k$. 
Let us consider the sequence $(\phi_{\sigma_n(n)})_{n\geq 0}$ and define $\hat\phi=\hat\phi_k$ on $(0,\pi-\delta_k)$. Thus $\phi_{\sigma_n(n)}\wkconv*{*}{}\hat\phi$ in $X_k$ for all $k\geq 0$. Moreover $\|\hat\phi\|_{X_k}\leq C$ for all $k\geq 0$ and so $\hat\phi:(0,\pi)\rightarrow\R$ is measurable and $\|\hat \phi\|_{\infty}\leq C$.
\end{demo}

\begin{rem}
The next example shows that the stronger conclusion $\phi_{n_k}\wkconv*{*}{}\hat\phi$ in $L^{\infty}(0,\pi)$ cannot be obtained. Let $\phi_n$ be the continuous function which vanishes on $[0,\pi-\frac{1}{n}]$ and grows linearly from $\pi-\frac{1}{n}$ to $\pi$ in such a way that $\int_0^{\pi}\phi_n=1$. By construction, for all $\delta>0$ there exists $N=N(\delta)$ such that $\|\phi_n\|_{L^{\infty}(0,\pi-\delta)}\leq 0$ if $n\geq N$. We prove that there is no subsequence $\phi_{n_k}$ converging in the weak$^*$ sense in $L^{\infty}(0,\pi)$. Otherwise there should exist $\hat\phi\in L^{\infty}(0,\pi)$ with $\phi_{n_k}\wkconv*{*}{}\hat\phi$. Testing this convergence with the constant function $\chi=1$ we deduce that
$\int_0^{\pi}\hat\phi = \lim_{k\rightarrow+\infty}\int_0^{\pi}\phi_{n_k}=1$. We can also use the characteristic function of the interval $[0,\pi-\frac{1}{N}]$ as the test function $\chi$. In this case we obtain $\int_0^{\pi-\frac{1}{N}}\hat\phi = \lim_{k\rightarrow+\infty}\int_0^{\pi-\frac{1}{N}}\phi_{n_k}=0$. Since $N$ is arbitrary these two facts are not compatible and so $\hat\phi$ cannot exist.
\end{rem}

\begin{lema}\label{lema:liminf-limsup}
Let $f:[0,+\infty)\rightarrow\R$ be a continuous function and let $(\phi_n)_{n\geq 0}$ be a sequence in $L^{\infty}(J)$ such that $\phi_n\rightarrow+\infty$ almost uniformly as $n\rightarrow+\infty$. This means that there exists a sequence of numbers $(M_n)$ with $M_n\rightarrow+\infty$ and such that $\phi_n(t)\geq M_n$ almost everywhere. In addition assume that, for some $\gamma\in L^{\infty}(J)$, $f\circ \phi_n\wkconv*{*}{}\gamma$ in $L^{\infty}(J)$. Then 
\begin{equation}\label{DV2}
\liminf_{x\rightarrow+\infty} f(x) \leq \gamma(t) \leq \limsup_{x\rightarrow+\infty}f(x)
\end{equation}
for almost every $t\in J$.
\end{lema}

\begin{demo}
If $\liminf_{x\rightarrow+\infty}f(x)$ is finite we take $\alpha\in\R$ with $\alpha<\liminf_{x\rightarrow+\infty}f(x)$. Otherwise we define $\alpha=-\infty$. In the same way we take an upper bound $\beta$ for $\limsup_{x\rightarrow+\infty}f(x)$. We define the set 
\[
\Omega_{\alpha,\beta}\!:=\{g\in L^{\infty}(J): \alpha\leq g(t) \leq \beta \text{ a.e.}\}.
\] 
The choice of $\alpha,\beta$ and the convergence to infinity of $\phi_n$ imply that $f\circ\phi_n$ belongs to $\Omega_{\alpha,\beta}$ if $n$ is large enough. The set $\Omega_{\alpha,\beta}$ is convex and closed with respect to the strong topology in $L^{\infty}(J)$. Thus, it is also closed with respect to the weak$^*$ topology. In consequence $\gamma$ belongs to $\Omega_{\alpha,\beta}$. The class of functions described by the inequalities~\eqref{DV2} is precisely the set $\bigcap_{\alpha,\beta}\Omega_{\alpha,\beta}$.
\end{demo}

\begin{prooftext}{Proof of Theorem~\ref{thm:limite-x}.}
Let us denote $f_I(t)\!:=\frac{1}{\sqrt{I}}x(t,I)$ and let $(I_n)_{n\geq 0}$ be a sequence of actions such that $I_n\rightarrow+\infty$ as $n\rightarrow+\infty$. The function $f_{I_n}$ satisfies the Cauchy problem
\begin{equation}\label{eq:Cauchy_f}
\begin{cases}
\ddot f_{I_n} + \dfrac{V'(f_{I_n}\sqrt{I_n})}{f_{I_n}\sqrt{I_n}}f_{I_n}=0,\\
f_{I_n}(0)=\frac{r_n}{\sqrt{I_n}},\ \dot f_{I_n}(0)=0,
\end{cases}
\end{equation}
where $r_n=x(0,I_n)$. Since $\Omega(I)=I$, according to $(a)$ and $(b)$ in Proposition~\ref{prop:cotas} we have that $\abs{f_{I_n}(t)}$ and $|\dot f_{I_n}(t)|\leq \sqrt{2}$ are uniformly bounded in $\T$ for all $n\geq 0$. Thus the family $\{f_{I_n}\}_{n\geq 0}$ is equicontinuous, bounded and $K$-Lipschitz with $K\leq\sqrt{2}$ so by Arzel\`a-Ascoli theorem there exists a subsequence $(I_{n_k})_{k\geq 0}\subset(I_n)_{n\geq 0}$ such that $f_{I_{n_k}}\rightarrow\psi$ uniformly for some function $\psi\in W^{1,\infty}(\T)$ with Lipschitz constant $[\psi]_{Lip}=K\leq\sqrt{2}$. This proves by now that $\psi\in W^{1,\infty}(\T)$.

By Lemma~\ref{lema:semiperiodo-cero}, we have that $\mathcal N(I_n)\rightarrow\{\pi\}$ as $n\rightarrow+\infty$. Consequently, for any $\delta>0$ small enough there exists $n_0(\delta)$ such that $[-\pi+\delta,\pi-\delta]\subset[-\pi,\pi]\setminus\mathcal N(I_n)$ for all $n\geq n_0(\delta)$. Thus, for each $\delta>0$, the function $f_{I_n}$ is positive on $[-\pi+\delta,\pi-\delta]$ for all $n\geq n_0(\delta)$ and so by hypothesis $(H_2)$ we have
\begin{equation}\label{eq:unif_bound}
\abs{\frac{V'(f_{I_n}(t)\sqrt{I_n})}{f_{I_n}(t)\sqrt{I_n}}}\leq \|V''\|_{\infty}
\ \text{ for all } t\in[-\pi+\delta,\pi-\delta] \text{ and }n\geq n_0(\delta).
\end{equation}
Thus, it follows from \eqref{eq:Cauchy_f} that there exists a constant $C>0$ such that $|\ddot f_{I_n}(t)|\leq C$ for all $t\in [-\pi+\delta,\pi-\delta]$, $\delta>0$ and $n\geq n_0(\delta)$. Again by Arzel\`a-Ascoli theorem, for each $\delta>0$ there exists a subsequence $(I_{n_k(\delta)})_{k\geq 0}$ such that $\dot f_{I_{n_k(\delta)}}\rightarrow\dot \psi$ uniformly on $[-\pi+\delta,\pi-\delta]$. Taking a sequence $(\delta_k)_{k\geq 0}$ tending to zero and arguing by a diagonal process we can consider the convergent subsequence independent from $\delta$ and so without loss of generality $\dot f_{I_{n_k}}\rightarrow\dot\psi$ uniformly in every compact subset of $(-\pi,\pi)$.

Up to now we have proved the convergence results $(a)$ and $(b)$. Moreover, we know by now that $\psi\in W^{1,\infty}(\T)$. Let us now show that $\psi$ is a solution of a bouncing problem. The first step will be the construction of $\gamma\in L^{\infty}(\T)$. On account of~\eqref{eq:unif_bound} and Lemma~\ref{lema:weak}, there exists $\gamma\in L^{\infty}(\T)$ with $\|\gamma\|_{L^{\infty}(\T)}\leq\|V''\|_{\infty}$ such that 
\begin{equation}\label{eq:weak-conv}
\frac{V'(f_{I_{n_k}}\sqrt{I_{n_k}})}{f_{I_{n_k}}\sqrt{I_{n_k}}}\wkconv*{\phantom{*}*\phantom{*}}{}\gamma\ \text{ in }L^{\infty}(-\pi+\delta,\pi-\delta)\text{ for any }\delta>0.
\end{equation}
For the sake of simplicity we consider $(I_{n_k})$ to be the same as the previous subsequence. From the previous limit we have that $\gamma$ is $2\pi$-periodic.

Let us take $\delta>0$ and $\phi\in\mathcal D(-\pi+\delta,\pi-\delta)$ a smooth test function with $\text{supp}(\phi)\subset[-\pi+\delta,\pi-\delta]$. We have that $\text{supp}(\phi)\subset\mathcal P(I_{n_k})$ for sufficiently large $k$. We have, on account of equation in \eqref{eq:Cauchy_f}, that
\[
\int_{-\pi}^{\pi}\left\{f_{I_{n_k}}(t)\ddot \phi(t)+\frac{V'(f_{I_{n_k}}(t)\sqrt{I_{n_k}})}{f_{I_{n_k}}(t)\sqrt{I_{n_k}}}f_{I_{n_k}}(t)\phi(t)\right\}dt=0.
\]
It is at this point where property \eqref{eq:weak-conv} plays an important role in the passage to the limit. Indeed, if $u_n\rightarrow u$ uniformly in $L^{\infty}(\T)$ and $v_n\wkconv*{*}{} v$ in $L^{\infty}(\T)$ with $\|v_n\|_{L^{\infty}(\T)}$ bounded then $\int_{\T} u_nv_n \rightarrow \int_{\T} uv$. Therefore, since $f_{I_{n_k}}$ converges uniformly to $\psi$ and using the weak$^*$ convergence in \eqref{eq:weak-conv},
\[
\int_{-\pi}^{\pi}\left\{f_{I_{n_k}}(t)\ddot \phi(t)+\frac{V'(f_{I_{n_k}}(t)\sqrt{I_{n_k}})}{f_{I_{n_k}}(t)\sqrt{I_{n_k}}}f_{I_{n_k}}(t)\phi(t)\right\}dt\rightarrow \int_{-\pi}^{\pi}
\{\psi(t)\ddot\phi(t) + \gamma(t)\psi(t)\phi(t)\}dt=0.
\]
That is, $\psi$ satisfies the ordinary differential equation $\ddot\psi+\gamma(t)\psi=0$ in the sense of distributions in $\mathcal D(-\pi,\pi)$. Since we are dealing with an o.d.e., distributional solutions coincide with solutions in the Carath\'eodory sense and so $\psi\in W^{2,\infty}(-\pi,\pi)$ and the equation holds almost everywhere in $(-\pi,\pi)$. This shows properties $(i)$ and $(ii)$ of Definition~\ref{defi:BP}.

Let us now show the remaining properties. To show $(iii)$ we select $t\in[0,2\pi]$, $t\neq \pi$. Then, since $\mathcal N(I_{n_k})\rightarrow\{\pi\}$, we know that $f_{I_{n_k}}(t)\geq 0$ for $k$ large enough. In consequence $\psi(t)=\lim_{k\rightarrow+\infty}f_{I_{n_k}}(t)$ is non-negative. Since $\psi$ is continuous and periodic we deduce that $\psi\geq 0$ in $\T$. Moreover, by Proposition~\ref{prop:cotas} we have $\abs{f_I(\pi)}\leq \abs{a} I^{-\frac{1}{2}}$. Thus $\abs{\psi(\pi)}\leq 0$ and so $\psi(\pi)=0$. On account of the $2\pi$-periodicity we have then that $\psi(t)=0$ if $t=(2n+1)\pi$. On the other hand, if there exists $t^*\in(-\pi,\pi)$ such that $\psi(t^*)=0$, since $\psi\geq 0$ we also have $\dot\psi(t^*)=0$. Notice that we already know that $\psi$ is in $W^{2,\infty}(-\pi,\pi)$ and so it also belongs to $C^1[-\pi,\pi]$. Thus, on account of uniqueness of solution of the differential equation $\ddot\psi+\gamma(t)\psi=0$ on $(-\pi,\pi)$ we have that $\psi\equiv 0$. However, by Lemma~\ref{lema:amplitud} we have $f_{I_n}(0)\geq \sqrt{\frac{2}{\|V''\|_{\infty}}}$ contradicting that $\psi\equiv 0$. This proves $(iii)$. The symmetry properties of the equation~\eqref{eq:system} lead to the identity $f_I(\pi+t)=f_I(\pi-t)$. Passing to the limit we obtain the identity $\psi(\pi+t)=\psi(\pi-t)$ which implies $\dot\psi(\pi+t)=-\dot\psi(\pi-t)$ for all $t\in[-\pi,\pi]$. Note that, by periodicity, $\psi$ also belongs to $C^1[\pi,3\pi]$. To complete the proof we must check the condition~\eqref{lab}. Let us fix a compact interval $J\subset(-\pi,\pi)$. Since $\psi(t)$ is positive on $J$ and $f_{I_n}\rightarrow\psi$ uniformly, there exists a constant $\mu_J>0$ such that $f_{I_n}(t)\geq \mu_J$ if $n$ is large enough and $t\in J$. We apply Lemma~\ref{lema:liminf-limsup} with $f(x)=\frac{V'(x)}{x}$ and $\phi_n(t)=f_{I_n}(t)\sqrt{I_n}$. Then $\phi_n\rightarrow+\infty$ almost uniformly on $J$ with $M_n=\mu_J\sqrt{I_n}$. Therefore~\eqref{lab} holds for almost every $t\in J$. The interval $J\subset(-\pi,\pi)$ is arbitrary and the inequality~\eqref{lab} will also hold almost everywhere on $(-\pi,\pi)$.
\end{prooftext}

\begin{prooftext}{Proof of Theorem~\ref{thm:limite-xI}.}
It is very similar to the proof of Theorem~\ref{thm:limite-x} and we will stress the points where there appear some differences. Let us denote $g_I(t)\!:=\sqrt{I}\frac{\partial x}{\partial I}(t,I)$ and let $(I_n)_{n\geq 0}$ be a sequence of actions such that $I_n\rightarrow+\infty$ as $n\rightarrow+\infty$. The function $g_{I_n}$ satisfies the Cauchy problem
\begin{equation}\label{eq:Cauchy-partialI}
\begin{cases}
\ddot g_{I_n} + V''(x(t,I_n)) g_{I_n} =0,\\
g_{I_n}(0)=\frac{\sqrt{I_n}}{V'(r_n)},\ \dot g_{I_n}(0)=0.\\
\end{cases}
\end{equation}
Since $\mathcal A(I_n)\subset\mathcal N(I_n)$, by Lemma~\ref{lema:semiperiodo-cero} we have $\mathcal A(I_n)\rightarrow \{\pi\}$ as $I_n\rightarrow+\infty$. Thus, for any $\delta>0$ small enough, there exists $n_0(\delta)$ such that $[-\pi+\delta,\pi-\delta]\subset [-\pi,\pi]\setminus\mathcal A(I_n)$ for all $n\geq n_0(\delta)$. Therefore, by definition of $\mathcal A(I)$, $\abs{V''(x(t,I_n))}\leq \|V''\|_{\infty}$ for all $t\in[-\pi+\delta,\pi-\delta]$ and, by $(c)$ in Proposition~\ref{prop:cotas}, $\abs{g_{I_n}(t)}$ is uniformly bounded in $\T$ for all $n\geq 0$. 
Consequently, taking into account equation \eqref{eq:Cauchy-partialI}, we have that there exists some constant $C>0$ such that $\|\ddot g_{I_n}\|_{L^{\infty}([-\pi+\delta,\pi-\delta])}\leq C$ and $\|\dot g_{I_n}\|_{L^{\infty}([-\pi+\delta,\pi-\delta])}\leq C$ for all $n\geq n_0(\delta)$. We point out here that $C$ is independent from $\delta$ since the previous bounds are independent from $\delta$. 
By Arzel\`a-Ascoli theorem, there exist a subsequence $(I_{n_k})_{k\geq 0}$ of $(I_n)_{n\geq n_0}$ and $\psi_1\in C^1(-\pi,\pi)$ such that $g_{I_{n_k}}\rightarrow\psi_1$ and $\dot g_{I_{n_k}}\rightarrow\dot\psi_1$ uniformly on compact subsets of $(-\pi,\pi)$. Moreover $|\psi_1(t)|\leq C$ and $|\dot\psi_1(t)|\leq C$ for all $t\in(-\pi,\pi)$. In particular $\psi_1$ is bounded and Lipschitz-continuous with $[\psi_1]_{\text{Lip}}\leq C$. This implies that $\psi_1$ has a continuous extension to $[-\pi,\pi]$ and, by periodicity, to the whole real line.

By now we have proved that $\psi_1\in W^{1,\infty}(\T)$ and the convergence assertion on the statement. The rest of the proof goes along the same lines of the proof of Theorem~\ref{thm:limite-x}. The only substantial difference is in the last step, when Lemma~\ref{lema:liminf-limsup} is applied. Now $f(x)=V''(x)$ and $\phi_n(t)=x(t,I_n)$. For each compact interval $J\subset(-\pi,\pi)$ the sequence $\phi_n\rightarrow+\infty$ uniformly in $J$.
\end{prooftext}

\section{Resonance of isochronous oscillators}\label{sec:resonance}

This section is dedicated to prove some of the results presented in Section~\ref{sec:results}. 

\subsection{Proof of Theorem~\ref{main}}\label{sec:main}

We first study some simple properties of the linear equation
\begin{equation}\label{lineq}
\ddot y + a(t) y = b(t)
\end{equation}
where $a\in L^{\infty}(0,2\pi)$, $b\in L^1(0,2\pi)$ and $\|a\|_{L^{\infty}(0,2\pi)}\leq A$.
\begin{enumerate}[$(i)$]
\item Let $\psi(t)$ be the solution of~\eqref{lineq} with $b\equiv 0$ and initial conditions $\psi(0)=1,\dot\psi(0)=i$. Then there exists $C>0$, depending only on $A$, such that $|\psi(t)|\leq C$ for each $t\in[0,2\pi]$.
\item Let $y(t)$ be the solution of~\eqref{lineq} with initial conditions $y(0)=\dot y(0)=0$. Then $|y(t)|\leq C^2\|b\|_{L^1(0,2\pi)}$ if $t\in[0,2\pi]$.
\end{enumerate}

To prove $(i)$ we observe that the matrix $\left(\begin{smallmatrix}u&v\\ \dot u&\dot v\end{smallmatrix}\right)$ is a solution of the linear system $\dot Y = \left(\begin{smallmatrix}0&1\\ -a(t)& 0\end{smallmatrix}\right) Y$, where $u=Re\psi$, $v=Im\psi$. The estimate is obtained as in Lemma~\ref{lema:bounds2}. To prove $(ii)$ we apply the formula of variation of constants, leading to
\[
y(t)=\int_0^t [v(t)u(s)-v(s)u(t)]b(s)ds.
\]
Cauchy-Schwarz inequality can be applied to the term in brackets to obtain
\[
|v(t)u(s)-v(s)u(t)|\leq (u(t)^2+v(t)^2)^{\frac{1}{2}}(u(s)^2+v(s)^2)^{\frac{1}{2}}\leq C^2.
\]

The Second Massera's Theorem~\cite{Massera} states that if all solutions of~\eqref{per} are globally defined in the future and at least one of them is bounded in the future, then a $2\pi$-periodic solution must exist. Our strategy will be to prove non-existence of $2\pi$-periodic solutions of~\eqref{per} for $\epsilon\neq 0$ small. We proceed by contradiction and assume that $\{\epsilon_n\}$ is a sequence with $\epsilon_n\neq 0$, $\epsilon_n\rightarrow 0$ such that the equation~\eqref{per} has a $2\pi$-periodic solution $x_n(t)$ for $\epsilon=\epsilon_n$. Let $X_n(t)$ be the solution of~\eqref{aut} with the same initial conditions of $x_n(t)$ at $t=0$, $X_n(0)=x_n(0)$ and $\dot X_n(0)=\dot x_n(0)$. The difference $y_n(t)=x_n(t)-X_n(t)$ can be seen as a solution of the equation~\eqref{lineq} with
\[
a(t)=\int_0^1 V''\bigl( (1-\lambda)x_n(t)+\lambda X_n(t) \bigr) d\lambda,\ b(t)=\epsilon_n p(t).
\]
From property $(ii)$ we deduce that
\begin{equation}\label{ebound}
\|y_n\|_{L^{\infty}(\T)}\leq C^2|\epsilon_n|\|p\|_{L^1(\T)}.
\end{equation}
Here $A=\sup_{x\in\R}|V''(x)|$.

The function $y_n(t)$ can also be interpreted as a $2\pi$-periodic solution of the linear periodic equation
\[
\ddot y + V''(X_n(t)) y = \epsilon_n p_n(t) - q_n(t)
\]
with
\[
q_n(t)=y_n(t)\int_0^1\bigl[ V''\bigl((1-\lambda)x_n(t)+\lambda X_n(t)\bigr) - V''\bigl(X_n(t)\bigr)\bigr] d\lambda.
\]
Note that the coefficients of this equation are $2\pi$-periodic because the oscillator defined by~\eqref{aut} is globally isochronous. All solutions of~\eqref{aut} are of the type $\varphi(t-\theta,r)$. In particular, $X_n(t)=\varphi(t-\theta_n,r_n)$ for some $\theta_n\in[0,2\pi]$ and $r_n\geq 0$. The function $\psi(t-\theta_n,r_n)$ is a $2\pi$-periodic solution of the homogeneous equation $\ddot y + V''\bigl(X_n(t)\bigr) y=0$. Fredholm alternative can be applied to deduce that
\[
\epsilon_n\int_0^{2\pi} p(t)\psi(t-\theta_n,r_n)dt - \int_0^{2\pi} q_n(t)\psi(t-\theta_n,r_n)dt =0.
\]

In other words,
\[
\Phi_p(\theta_n,r_n)=\frac{1}{2\pi}\int_0^{2\pi}\frac{q_n(t)}{\epsilon_n}\psi(t-\theta_n,r_n)dt.
\]

The definition of $q_n$ together with the estimate~\eqref{ebound} and the Dominated Convergence Theorem imply that $\frac{1}{\epsilon_n} q_n(t)\rightarrow 0$ as $n\rightarrow+\infty$ for every $t\in[0,2\pi]$. From the property $(i)$ we know that $\|\psi(\cdot, r_n)\|_{L^{\infty}(\T)}\leq C$ and the estimate~\eqref{ebound} implies that
\[
\frac{\|q_n\|_{L^{\infty}(\T)}}{|\epsilon_n|}\leq 2C^2\|p\|_{L^1(\T)}\|V''\|_{L^{\infty}(\R)}.
\] 
Using once again dominated convergence,
\[
\Phi_p(\theta_n,r_n)\rightarrow 0 \text{ as }n\rightarrow+\infty.
\]
This is incompatible with the assumption~\eqref{rc} and proves~Theorem~\ref{main}.\Qed

A refinement of Massera's Theorem shows that the equation~\eqref{per} has a $2\pi$-periodic solution if there exists a solution $x(t)$ which is defined in $[0,+\infty)$ and satisfies
\[
\liminf_{n\rightarrow+\infty}\left(\abs{x(2\pi n)}+\abs{\dot x(2\pi n)}\right)<+\infty.
\]
See~\cite{Ortega3}. In consequence, if the equation~\eqref{per} is resonant all solutions satisfy
\[
\lim_{n\rightarrow+\infty}\left(\abs{x(2\pi n)}+\abs{\dot x(2\pi n)}\right)=+\infty.
\]
A standard argument using Gronwall's lemma and $\|V''\|_{\infty}<+\infty$ implies that there exists a constant $C>0$ such that for every solution $x(t)$,
\[
\abs{x(t)}+\abs{\dot x(t)}\geq C\left(\abs{x(0)}+\abs{\dot x(0)}\right)\ \text{ if }t\in[0,2\pi].
\]
Therefore every solution of~\eqref{per} satisfies $\abs{x(t)}+\abs{\dot x(t)}\rightarrow+\infty$ as $t\rightarrow+\infty$.

\subsection{Proof of Proposition~\ref{reciprocal}}\label{sec:reciprocal}

The aim is to show that if $\Phi_p$ has a non-degenerate zero $(\theta_*,r_*)$ then a $2\pi$-periodic solution persists on the perturbed system~\eqref{per} for small $\epsilon$.
At this point the authors want to emphasize that Melnikov method to study subharmonic bifurcations cannot be applied in this framework due to the isochronicity of the center. Instead, degree theory will be the key point on the proof of the result. For the sake of completeness we sketch Nagumo's definition of Brouwer degree. There are other equivalent ways to define this degree, we refer the reader to \cite{Lloyd} for more details.

Let $\Omega$ be an open and bounded subset of $\R^d$ and a function $f:\bar\Omega\rightarrow\R^d$ which is continuous and does not vanish on the boundary of $\Omega$. First we assume that $f\in C^1(\bar\Omega)$ and has a finite number of zeros $x_1,\dots, x_n\in\Omega$ with $\det f'(x_i)\neq 0$ for each $i$. Then the Brouwer degree of $f$ in $\Omega$ is
\[
\deg(f,\Omega)=\sum_{i=1}^n \text{sign}\det f'(x_i).
\]
In the case $f$ is continuous, we approximate it by functions $f_k$ in the previous conditions and define
\[
\deg(f,\Omega)=\lim_{k\rightarrow+\infty}\deg(f_k,\Omega).
\]

An important property of the degree is its invariance by homotopy, which plays a crucial role on the proof of the result under consideration.

\begin{prooftext}{Proof of Proposition~\ref{reciprocal}.}
Let $x(t;x_0,v_0,\epsilon)$ be the solution of~\eqref{per} with initial conditions $x(0)=x_0$, $\dot x(0)=v_0$. Define
\[
\Delta(t;x_0,v_0,\epsilon)\!:=\begin{cases}
\frac{1}{\epsilon}\bigl[x(t;x_0,v_0,\epsilon)-x(t;x_0,v_0,0)\bigr] & \text{ if }\epsilon\neq 0,\\
\frac{\partial x}{\partial \epsilon}(t;x_0,v_0,0) & \text{ if }\epsilon=0.
\end{cases}
\]
For each $\epsilon\in\R$ consider the planar map
\[
\mathcal F_{\epsilon}:\R^2\rightarrow\R^2,\ \mathcal F_{\epsilon}(x_0,v_0)=\bigl(\Delta(2\pi;x_0,v_0,\epsilon), \dot \Delta(2\pi;x_0,v_0,\epsilon)\bigr).
\]
We are going to prove that $\mathcal F_{0}$ has a non-degenerate zero. Then, for small $\epsilon\neq 0$, the map $\mathcal F_{\epsilon}$ will also have a zero. Clearly this zero will produce a $2\pi$-periodic solution of~\eqref{per}.

Let $N\geq 1$ be the integer such that all non-trivial solution of~\eqref{aut} have minimal period $\frac{2\pi}{N}$. Then $\psi(\cdot,r)$ has period $\frac{2\pi}{N}$ and the map
\[
\tilde\Phi_p(\theta,r)\!:=\Phi_p\bigl(\tfrac{\theta}{N},r\bigr)
\]
is well defined on the cylinder $\mathcal C$. The map $(\theta,r)\mapsto (N\theta,r)$ is a local diffeomorphism of the open cylinder $\dot{\mathcal C}=(\R/2\pi\Z)\times (0,\infty)$. Let $(\theta_*,r_*)$ be the non-degenerate zero of $\Phi_p$ given by assumption, then $(N\theta_*,r_*)$ is a non-degenerate zero of $\tilde\Phi_p$. This can be deduced from the definition of degree or from Leray multiplication theorem.

Let us now consider the diffeomorphism
\[
\eta:(\theta,r)\in\dot{\mathcal C}\mapsto\bigl(\varphi\bigl(-\tfrac{\theta}{N},r\bigr),\dot\varphi\bigl(-\tfrac{\theta}{N},r\bigr)\bigr)\in\R^2\setminus\{0\}.
\]
We are going to compute $\mathcal F_{0}\circ\eta$ using the formula of variation of constants. The function $y(t)=\Delta(t;\eta(\theta,r),0)$ is the solution of
\[
\ddot y + V''\bigl(\varphi\bigl(t-\tfrac{\theta}{N},r\bigr)\bigr) y = p(t),\ y(0)=\dot y(0)=0.
\]
Then
\[
\left(\begin{matrix}
\Delta(t;\eta(\theta,r),0)\\
\dot\Delta(t;\eta(\theta,r),0)
\end{matrix}\right)= M\bigl(t-\tfrac{\theta}{N},r\bigr)\int_0^t
\left(\begin{matrix}
u\bigl(s-\tfrac{\theta}{N}\bigr)\\
v\bigl(s-\tfrac{\theta}{N}\bigr)
\end{matrix}\right) p(s) ds
\]
where $M=\left(\begin{smallmatrix} v & -u\\ \dot v & -\dot u\\ \end{smallmatrix}\right)$ and $u(t,r)=Re\psi(t,r)$, $v(t,r)=Im\psi(t,r)$. In consequence
\[
\mathcal F_0 \circ \eta(\theta,r)=2\pi M\bigl(-\tfrac{\theta}{N},r\bigr)\tilde\Phi_p(\theta,r).
\]
Here we are identifying $\C$ and $\R^2$ so that $\tilde\Phi_p$ takes values in $\R^2$. The matrix $M$ has determinant one and so the zeros of $\mathcal F_0\circ\eta$ and $\tilde\Phi_p$ coincide. For each $\lambda\in[0,1]$ define $H_{\lambda}:\dot{\mathcal C}\rightarrow \R^2$ by $H_{\lambda}(\theta,r)=M(-\tfrac{\lambda\theta}{N},r)\tilde\Phi_p(\theta,r)$. We observe that $H_0=J \tilde\Phi_p$ and $H_1=\tfrac{1}{2\pi}\mathcal F_0\circ\eta$ where $J=\left(\begin{smallmatrix} 0 & -1\\ 1 & 0\\ \end{smallmatrix}\right)$. Moreover, the zeros of $H_{\lambda}$ are independent of $\lambda$. We conclude that $(N\theta_*,r_*)$ is a non-degenerate zero of $\mathcal F_0\circ\eta$. Then $(x_0^*,v_0^*)=\eta(N\theta_*,r_*)$ is a non-degenerate zero of $\mathcal F_0$. This completes the proof.
\end{prooftext}

\subsection{Proof of Proposition~\ref{functional}}\label{sec:functional}

\begin{lema}\label{lema:cotas_psi}
There exist $C,c>0$ depending on $\|V''\|_{\infty}$ such that
\[
c\leq \abs{\psi(t,r)}\leq C, \ |\dot\psi(t,r)|\leq C\ \text{ if }t\in[0,2\pi] \text{ and } r>0.
\]
\end{lema}
\begin{demo}
Assume that $\psi=u+iv$ with $u=\text{Re}\psi$, $v=\text{Im}\psi$. Then $\left(\begin{smallmatrix} u & v\\ \dot u & \dot v\\ \end{smallmatrix}\right)$ is a matrix solution of the linearized system and the upper estimates are obtained as in Lemma~\ref{lema:bounds2}. In particular $\dot u(t)^2+\dot v(t)^2\leq C^2$ if $t\in[0,2\pi]$. From Liouville's identity,
\[
1=u\dot v-\dot u v\leq (u^2+v^2)^{\frac{1}{2}}(\dot u^2+\dot v^2)^{\frac{1}{2}}\leq C\abs{\psi}
\]
and we can take $c=\frac{1}{C}$.
\end{demo}

\begin{prooftext}{Proof of Proposition~\ref{functional}}
Let $p,q\in L^1(\T)$. The fact that $\mathcal R_V$ is open follows from Lemma~\ref{lema:cotas_psi} since
\[
\abs{\Phi_p(\theta,r)-\Phi_q(\theta,r)}\leq \frac{C}{2\pi}\|p-q\|_{L^1(\T)}.
\]
Let us now show that $\mathcal R_V$ is non-empty. To do so, we shall see that if we take a sequence $p_n\in L^1(\T)$ with $p_n\wkconv*{}{} \delta$ then $p_n\in\mathcal R_V$ for $n$ large enough. In order to arrive at contradiction we consider a sequence $\{p_n\}_{n\geq 0}$ such that for all $n\geq 0$, $p_n\in L^1(\T)$, $p_n\notin\mathcal R_V$ and $p_n\wkconv*{}{} \delta$. From $p_n\notin\mathcal R_V$ it follows that for each $n$ there exists a sequence $(\theta_k^{(n)},r_k^{(n)})\in\mathcal C$ such that $\Phi_{p_n}(\theta_k^{(n)},r_k^{(n)})\rightarrow 0$. Arguing by a diagonal process we can consider  $\Phi_{p_n}(\theta_n,r_n)\rightarrow 0$. On account of the bounds in Lemma~\ref{lema:cotas_psi}, by Arzel\`a-Ascoli theorem, there exists a subsequence $\psi(t+\theta_{n_k},r_{n_k})$ that converges uniformly on $[0,2\pi]$ to some $\hat\psi\in C(\T,\C)$ as $k\rightarrow+\infty$. Moreover, $c\leq|\hat\psi(t)|\leq C$. Thus,
\[
2\pi\Phi_{p_{n_k}}(\theta_{n_k},r_{n_k})=\int_0^{2\pi}p_{n_k}(t)\psi(t+\theta_{n_k},r_{n_k})dt=
\int_0^{2\pi}p_{n_k}(t)\left(\psi(t+\theta_{n_k},r_{n_k})-\hat\psi(t)\right)dt+\int_0^{2\pi}p_{n_k}(t)\hat\psi(t)dt.
\]
The first integral tends to zero as $k$ tends to infinity. Indeed, we have
\[
\int_0^{2\pi}p_{n_k}(t)\left(\psi(t+\theta_{n_k},r_{n_k})-\hat\psi(t)\right)dt\leq \|p_{n_k}\|_{L^1(\T)}\|\psi(\cdot+\theta_{n_k},r_{n_k})-\hat\psi\|_{L^{\infty}(\T)}\rightarrow 0
\]
as $k\rightarrow+\infty$ due to the boundedness of $\|p_{n_k}\|_{L^1(\T)}$ and the convergence of $\psi(t+\theta_{n_k},r_{n_k})$ to $\hat\psi(t)$ uniformly in $[0,2\pi]$. In addition, since $p_n\wkconv*{}{}\delta$,
\[
\int_0^{2\pi}p_{n_k}(t)\hat\psi(t)dt\rightarrow \hat\psi(0)
\]
as $k\rightarrow+\infty$ with $|\hat\psi(0)|\geq c>0$. Consequently, 
\[
|\Phi_{p_{n_k}}(\theta_{n_k},r_{n_k})|\rightarrow \frac{1}{2\pi}|\hat\psi(0)|\geq \frac{c}{2\pi}>0
\]
as $k\rightarrow+\infty$, reaching contradiction with the choice of $p_n\notin \mathcal R_V$.
\end{prooftext}

\section{Two resonance results for a periodic perturbed isochronous of Pinney type}\label{sec:pinney}

This section is dedicated to prove Theorem~\ref{thm:resonancia_pinney} and~\ref{thm:multiplicativo}.

\subsection{Proof of Theorem~\ref{thm:resonancia_pinney}}\label{aditiva}

We are concerned with the potential
\begin{equation}\label{pinn}
V(x)=\frac{1}{8}\left((x+1)^2+\frac{1}{(x+1)^2}\right)-\frac{1}{4}, \ x\in\mathcal J=(-1,+\infty).
\end{equation}
All the assumptions $(H_i)$ are satisfied in this case. From the first derivative $V'(x)=\frac{1}{4}\left(x+1-\frac{1}{(x+1)^3}\right)$ it is easy to deduce that $(H_0)$ and $(H_3)$ hold. Note that $\frac{V'(x)}{x}\rightarrow\frac{1}{4}$ as $x\rightarrow+\infty$. The second derivative $V''(x)=\frac{1}{4}+\frac{3}{4}(x+1)^{-4}$ is positive and so $(H_1)$ also holds. Finally we observe that $V''(x)\leq 1$ if $x\geq 0$ and this implies $(H_2)$.

The next result could be stated for a general class of isochronous systems with strictly convex potential. For our purposes we can restrict to the potential defined by~\eqref{pinn}.

\begin{lema}\label{lema:partial_I_decrece}
In the previous notations assume that $x(t,I)$ is a solution of~\eqref{pin}. Then
\[
\frac{\partial \dot x}{\partial I}(t,I)<0 \text{ on }(0,\pi).
\]
\end{lema}
\begin{demo}
The functions $\dot x(t,I)$ and $\frac{\partial x}{\partial I}(t,I)$ are a fundamental pair of solutions of the variational equation $\ddot\xi+V''(x(t,I))\xi=0$. From equation~\eqref{aut} we deduce that the critical points of $\dot x(t,I)$ coincide with the zeros of $x(t,I)$. In consequence $\dot x(t,I)$ has exactly two critical points in the interval $[-\pi,\pi]$, say $-\pi<t^*_-<0<t^*_+<\pi$. Moreover, by symmetry, $t^*_-=-t^*_+$. We are going to apply Lemma~\ref{lema:2_criticos} twice. First we take $a(t)=V''(x(t,I))$, $\phi_1(t)=\frac{\partial x}{\partial I}(t,I)$, $\phi_2(t)=\dot x(t,I)$, $\tau=t^*_+$. Since $\dot x(t,I)$ has no critical points on $(t^*_-,t^*_+)$ we deduce that $\frac{\partial x}{\partial I}(t,I)$ has no critical points on $(t^*_-,t^*_+)\setminus\{0\}$. Next we observe that $x(t,I)$ has the symmetry, $x(t+\pi,I)=x(-t+\pi,I)$. This allows us to apply again the Lemma with $\tilde a(t)=a(t+\pi)$, $\tilde\phi_1(t)=\phi_1(t+\pi)$, $\tilde\phi_2(t)=\phi_2(t+\pi)$, $\tau=\pi-t^*_+$. Now we conclude that $\frac{\partial x}{\partial I}(t,I)$ has no critical points on $(t^*_+,2\pi+t^*_-)\setminus\{\pi\}$. We also observe that $t^*_+$ cannot be a critical point of $\frac{\partial x}{\partial I}(t,I)$ because this function and $\dot x(t,I)$ are linearly independent. Summing up the previous discussion we conclude that $\frac{\partial x}{\partial I}(t,I)$ has no critical points on the interval $(0,\pi)$. In consequence $\frac{\partial x}{\partial I}(t,I)$ is monotone on this interval. Finally, by Lemma~\ref{lema:x_partial_I}, $\frac{\partial x}{\partial I}(0,I)>0$ and $\frac{\partial x}{\partial I}(\pi,I)<0$ and so $\frac{\partial\dot x}{\partial I}(t,I)$ is decreasing. The conclusion follows.
\end{demo}

The next result describes the asymptotic behaviour of $x(t,I)$ and its derivatives. It was already obtained in~\cite{BFS} but there are significant differences in the type of convergence and also in the proof.

\begin{prop}\label{prop:limites-uniformes}
Let $x(\theta,I)$ be a solution of system \eqref{pin}. Then, as $I\rightarrow+\infty$,
\begin{enumerate}[$(a)$]
\item $\frac{1}{\sqrt{I}}\ x(t,I)\rightarrow 2\sqrt{2}\abs{\cos\bigl(\frac{t}{2}\bigr)}$ uniformly in $\T$.
\item $\frac{1}{\sqrt{I}}\ \dot x(t,I)\rightarrow -\sqrt{2}\text{sign}(\cos\bigl(\frac{t}{2}\bigr))\sin\bigl(\frac{t}{2}\bigr)$ uniformly in every compact subset of $(-\pi,\pi)$.
\item $\sqrt{I}\ \frac{\partial x}{\partial I}(t,I)\rightarrow \sqrt{2}\abs{\cos\bigl(\frac{t}{2}\bigr)}$ uniformly in $\T$.
\end{enumerate}
\end{prop}

\begin{demo}
We apply the results of Section~\ref{sec:comportamiento}. A key observation is that the inequality in~\eqref{lab} now becomes an identity because the limit exists. Namely, $\lim_{x\rightarrow+\infty}\frac{V'(x)}{x}=\frac{1}{4}$. In consequence the function $\gamma$ has to be the constant $\frac{1}{4}$. The next step will be to determine all the solutions of the Bouncing problem when $\gamma=\frac{1}{4}$. A direct computation shows that these solutions are of the type
\[
\psi(t)=\lambda\abs{\cos\left(\frac{t}{2}\right)} \text{ with }\lambda>0.
\]
To determine the admissible values of $\lambda$ we first compute a couple of limits. Recalling that $x(0,I)=r$ with $V(r)=I$, $r>0$, and using the concrete form of the potential given by~\eqref{pinn} we deduce that 
\[
\lim_{I\rightarrow+\infty}\frac{1}{\sqrt{I}}x(0,I)=2\sqrt{2}.
\]
Also, from the identity $\frac{\partial x}{\partial I}(0,I)=\frac{1}{V'(r)}$, 
\begin{equation}\label{l2}
\lim_{I\rightarrow+\infty}\sqrt{I}\frac{\partial x}{\partial I}(0,I)=\sqrt{2}.
\end{equation}
The existence of these limits implies that the functions $\psi$ and $\psi_1$ are uniquely determined. Actually $\psi(t)=2\sqrt{2}\abs{\cos\left(\frac{t}{2}\right)}$ and $\psi_1(t)=\frac{1}{2}\psi(t)$. In consequence, as $I\rightarrow+\infty$, there is convergence of the functions $\frac{1}{\sqrt{I}}x(\cdot,I)$, $\frac{1}{\sqrt{I}}\dot x(\cdot,I)$ and $\sqrt{I}\frac{\partial x}{\partial I}(\cdot,I)$. The convergence must be understood in the senses specified by Theorems~\ref{thm:limite-x} and~\ref{thm:limite-xI}.

Assertions $(a)$ and $(b)$ are already proved but some additional work is required to prove $(c)$. In view of Theorem~\ref{thm:limite-xI} we know that $\sqrt{I}\frac{\partial x}{\partial I}(t,I)$ converges to $\psi_1(t)$ for each $t\in(-\pi,\pi)$. Moreover, by item $(b)$ of Lemma~\ref{lema:cota_xI_pi},
\[
\lim_{I\rightarrow+\infty}\sqrt{I}\frac{\partial x}{\partial I}(\pi,I)=\psi_1(\pi)=0.
\]
By symmetry we know that there is pointwise convergence of $\sqrt{I}\frac{\partial x}{\partial I}(t,I)$ to $\psi_1(t)$ for each $t\in\R$.  We must prove that this convergence is uniform and it is enough to consider the interval $[0,\pi]$. We know by Lemma~\ref{lema:partial_I_decrece} that the function $\sqrt{I}\frac{\partial x}{\partial I}(\cdot,I)$ is monotone in $[0,\pi]$. A classical result says that if a continuous function $\psi_1(t)$ is the pointwise limit of a sequence of monotone functions $\sqrt{I}\frac{\partial x}{\partial I}(\cdot,I)$, then the convergence is also uniform. 
\end{demo}

This result can be proved in a more direct way using the explicit formulas for $\varphi(t,r)$ and $\psi(t,r)$ given in Section~\ref{sec:results}. The proof above is more flexible and can be adapted to general families of potentials having~\eqref{pinn} as a prototype. 

We finish these preliminary results with some estimates for $x(t,I)$ and its derivatives.

\begin{prop}\label{prop:cotas-isocrono}
Assume that $V$ is given by~\eqref{pinn} and let $x(t,I)$ be the solution corresponding to~\eqref{aut}. Then there exists a constant $C>0$ independent of $I$ such that
\[
\frac{1}{\sqrt{I}}\abs{x(t,I)}+\frac{1}{\sqrt{I}}\abs{\dot x(t,I)} + \sqrt{I}\abs{\frac{\partial x}{\partial I}(t,I)}\leq C
\]
it $t\in[0,2\pi]$ and $I>0$.
\end{prop}
\begin{demo}
The bounds on $\abs{x(t,I)}$ and $\abs{\dot x(t,I)}$ are direct consequences of items $(a)$ and $(b)$ in Proposition~\ref{prop:cotas}. To get the bound on $\abs{\frac{\partial x}{\partial I}(t,I)}$ we can combine item $(c)$ with the continuity of the function
\[
(t,I)\in[0,2\pi]\times[0,+\infty)\mapsto \sqrt{I}\frac{\partial x}{\partial I}(t,I)\in\R.
\]
Note that this continuity is a consequence of the theorem of continuous dependence with respect to parameters and initial conditions when it is applied to the equation $\ddot y+V''(x(t,I))y=0$. The continuity at $I=0$ is also valid because $\lim_{I\rightarrow 0^+}x(t,I)=0$ uniformly in $t$ and $\lim_{I\rightarrow 0^+}\sqrt{I}\frac{\partial x}{\partial I}(0,I)=\sqrt{2}.$
\end{demo}

\begin{prooftext}{Proof of Theorem~\ref{thm:resonancia_pinney}}
We write the equation~\eqref{per} as the first order system $\dot x=y$, $\dot y = -V'(x)+\epsilon p(t)$. This system is defined on $D=\{(x,y)\in\R^2: x>-1\}$ and all solutions are global. In fact we can repeat the energy arguments of Remark~\ref{rem:def-global} for the potential $V(x)$ given by~\eqref{pinn}. The estimate~\eqref{EEst} now implies that, in finite time, the norm of the solutions $\abs{x(t)}+\abs{y(t)}$ cannot blow up and also that $x(t)$ cannot approach the vertical line $x=-1$. Therefore the maximal solutions $(x(t),y(t))$ are defined for all $t\in\R$.

Second Massera's Theorem is not directly applicable because the system is not defined in the whole plane. To overcome this difficulty we transport the system from $D$ to the plane $(z,y)$ with $x+1=e^z$. In the new variables the system is $\dot z=e^{-z}y$, $\dot y = -V'(e^z-1)+\epsilon p(t)$ and Massera's Theorem is applicable. We can conclude that if the equation~\eqref{per} has no $2\pi$-periodic solutions then every solution $x(t)$ of \eqref{per} satisfies
\begin{equation}\label{psb}
\limsup_{t\rightarrow+\infty}\left(\abs{x(t)}+\abs{\dot x(t)}+\frac{1}{x(t)+1}\right)=+\infty.
\end{equation}
We claim that every solution of~\eqref{per} satisfying~\eqref{psb} has to be unbounded, for otherwise there should exist a number $C>0$, a sequence $t_n\rightarrow+\infty$ and a solution $x(t)$ with $\abs{x(t)}+\abs{\dot x(t)}\leq C$ if $t\geq 0$ and $x(t_n)\rightarrow -1$. Then, for each $h\in[0,1]$, $1+x(t_n+h)\leq 1+x(t_n)+Ch$ and so
\[
\int_{t_n}^{t_n+1}\frac{dt}{(x(t)+1)^3}\geq \int_0^1\frac{dh}{(x(t_n)+Ch+1)^3}\rightarrow+\infty.
\]
After integrating the equation~\eqref{per} over the interval $[t_n,t_n+1]$ we obtain
\[
\dot x(t_n+1)-\dot x(t_n)+\frac{1}{4}\int_{t_n}^{t_n+1}(x(t)+1)dt=\epsilon\int_{t_n}^{t_n+1}p(t)dt+ \frac{1}{4}\int_{t_n}^{t_n+1}\frac{dt}{(x(t)+1)^3}.
\]
This identity cannot hold for large $n$ because the last term goes to infinity while the rest remain bounded.

The previous discussions allow us to conclude that the equation~\eqref{per} will be resonant as soon as it has no $2\pi$-periodic solutions. Our strategy will be to prove non-existence of $2\pi$-periodic solutions of~\eqref{per} for $\epsilon\neq 0$ small. We will proceed in several steps. First we prove that the initial condition of a $2\pi$-periodic solution cannot be close to the origin.

\begin{claim}
There exist $\sigma>0$ and $\epsilon_*>0$ such that if $0<\abs{\epsilon}<\epsilon_*$ then every $2\pi$-periodic solution of~\eqref{per} satisfies $\abs{x(0)}+\abs{\dot x(0)}\geq \sigma$.
\end{claim}

To prove this claim we denote by $x(t;x_0,v_0,\epsilon)$ the solution of~\eqref{per} with initial conditions $x(0)=x_0$, $\dot x(0)=v_0$ and define
\begin{equation}\label{eqdeg}
\mathcal F_{\epsilon}(x_0,v_0)\!:=\left(
\begin{matrix}
x(2\pi;x_0,v_0,\epsilon)-x_0\\
\dot x(2\pi;x_0,v_0,\epsilon)-v_0\\
\end{matrix}\right).
\end{equation}
The initial condition of a $2\pi$-periodic solution corresponds to a zero of $\mathcal F_{\epsilon}$. After differentiating with respect to the parameter $\epsilon$ we obtain the expansion
\begin{equation}\label{ae1}
\mathcal F_{\epsilon}(x_0,v_0)=\epsilon Y(2\pi;x_0,v_0)+o(\epsilon), \ \epsilon\rightarrow 0
\end{equation}
where $Y(t;x_0,v_0)$ is the solution of the linear equation
\[
\dot Y = \left(\begin{matrix} 0 & 1\\ -V''(x(t;x_0,v_0)) & 0\\ \end{matrix}\right) Y + \left(\begin{matrix} 0\\ p(t)\\ \end{matrix}\right),\ Y(0)=0.
\]
Moreover this expansion is uniform on $(x_0,v_0)\in K$ if $K$ is a compact subset of $D$. For the origin $(x_0,v_0)=(0,0)$ this system can be solved explicitly and a simple computation shows that
\[
Y(2\pi;0,0)=\int_0^{2\pi} p(t)\left(\begin{matrix} -\sin(t)\\ \cos(t) \\ \end{matrix}\right)dt.
\]
The condition~\eqref{rc} at $r=0$ implies that this vector does not vanish and, by continuous dependence, we can find $\sigma>0$ such that
\[
\|Y(2\pi;x_0,v_0)\|\geq \frac{1}{2}\|Y(2\pi;0,0)\|>0
\]
if $\abs{x_0}+\abs{v_0}\leq \sigma$. The expansion~\eqref{ae1} allows to find $\epsilon_*>0$ such that $\mathcal F_{\epsilon}(x_0,v_0)\neq 0$ if $0<\abs{\epsilon}<\epsilon_*$ and $\abs{x_0}+\abs{v_0}\leq \sigma$, proving the veracity of the claim.

Once we have proved the above claim we will employ action-angle variables $(\theta,I)$. In the punctured plane $\R^2\setminus\{0\}$ the Hamiltonian system~\eqref{per} is transformed into
\begin{equation}\label{eq:aa-var}
\begin{cases}
\dot\theta = 1+\epsilon p(t)\dfrac{\partial x}{\partial I}(\theta,I),\\
\dot I = -\epsilon p(t)\dfrac{\partial x}{\partial \theta}(\theta,I),
\end{cases}
\end{equation}
with Hamiltonian function $\mathscr H(\theta,I,t)=I+\epsilon p(t)x(\theta,I)$, where $x(t,I)$ is the solution of~\eqref{eq:system} introduced in Section~\ref{sec:4}.

In view of the Claim 1 every $2\pi$-periodic solution of~\eqref{per} corresponds to a solution $(\theta(t),I(t))$ of~\eqref{eq:aa-var} but the periodicity of this new solution has to be proved.

\begin{claim}
There exists $\epsilon_{**}>0$ such that if $0<\abs{\epsilon}<\epsilon_{**}$ and $x(t)$ is a $2\pi$-periodic solution of~\eqref{per} then $(\theta(t),I(t))$ is well defined for each $t\in\R$ and satisfies
$
\theta(t+2\pi)=\theta(t)+2\pi,\ I(t+2\pi)=I(t).
$
\end{claim}

Let $[0,\omega)$ be the maximal interval of $(\theta(t),I(t))$ when it is understood as a solution of~\eqref{eq:aa-var}. In view of the second equation~\eqref{eq:aa-var} and Proposition~\ref{prop:cotas-isocrono} we deduce that, for each $t\in[0,\omega)$,
\begin{equation}
\abs{\frac{d}{dt}(\sqrt{I(t)})}\leq \frac{C}{2}\abs{\epsilon}\abs{p(t)}.
\end{equation}
Integrating this inequality and assuming that $t\in[0,\omega)\cap[0,2\pi]$,
\begin{equation}\label{l111}
\abs{\sqrt{I(t)}-\sqrt{I(0)}}\leq \frac{C}{2}\abs{\epsilon} \|p\|_{L^1(\T)}.
\end{equation}
This estimate shows that $I(t)$ cannot blow up on the interval $[0,2\pi]$. Let us show also that $I(t)$ cannot touch the singularity $I=0$. For this purpose we select a closed orbit $\gamma$ of~\eqref{eq:system} contained in $|x|+|\dot x|<\sigma$ and denote by $\sigma_*>0$ the area enclosed by $\gamma$. Assuming that $\frac{C}{2}|\epsilon|\|p\|_{L^1(\T)}<\frac{1}{2}\sigma_*$ we deduce from Claim $1$ and the condition~\eqref{l111} that
\begin{equation}\label{bnec}
I(t)\geq \frac{1}{2}\sigma_*.
\end{equation}
In principle this estimate is valid for $t\in[0,\omega)\cap[0,2\pi]$ but now it is clear that $\omega>2\pi$. Since $(\theta(t),I(t))$ comes from a $2\pi$-periodic solution of~\eqref{per} we know that
\[
\theta(2\pi)=\theta(0)+2\pi M,\ I(2\pi)=I(0)
\]
for some integer $M$. The system~\eqref{eq:aa-var} is $2\pi$-periodic in $t$ and $\theta$ and we can conclude that $(\theta(t),I(t))$ can be extended to a global solution satisfying
\[
\theta(t+2\pi)=\theta(t)+2\pi M,\ I(t+2\pi)=I(t),\ t\in\R.
\]
To complete the proof of the claim we must show that $M=1$. From the first equation in~\eqref{eq:aa-var}, Proposition~\ref{prop:cotas-isocrono} and~\eqref{bnec} we deduce
\begin{equation}\label{EST2}
\abs{\frac{d}{dt}(\theta(t)-t)}\leq \frac{C|\epsilon|\abs{p(t)}}{\sqrt{I(t)}}\leq \sqrt{\frac{2}{\sigma_*}}C|\epsilon|\abs{p(t)},\ t\in[0,2\pi].
\end{equation}
For small $\epsilon$ the function $\theta(t)-t-\theta(0)$ must be small on $[0,2\pi]$ and this implies that $M=1$.

\begin{claim}
For small $\epsilon\neq 0$ the equation~\eqref{per} has no $2\pi$-periodic solutions.
\end{claim}

By an indirect argument we assume that $\{\epsilon_n\}_{n\geq 0}$, $\epsilon_n\neq 0$, is a sequence tending to zero such that $x_n(t)$ is a $2\pi$-periodic solution of~\eqref{per} with $\epsilon=\epsilon_n$. Let $(\theta_n(t),I_n(t))$ be the corresponding solution of~\eqref{eq:aa-var}. After extracting a subsequence we can assume
$I_n(0)\rightarrow I_*\in[0,+\infty]$ and $\theta_n(0)\rightarrow \theta_*\in[0,2\pi]$. In view of Claim $2$ we know that $\theta_n(t+2\pi)=\theta_n(t)+2\pi$, $I_n(t+2\pi)=I_n(t)$.

Let us assume first that $I_*<+\infty$. From Claim $1$ we know that $I_*\geq \sigma_*$ and the estimates~\eqref{l111} and~\eqref{EST2} imply that $I_n(t)\rightarrow I_*$ and $\theta_n(t)-t\rightarrow\theta_*$ uniformly in $t\in[0,2\pi]$. In addition, on account of equations in system \eqref{eq:aa-var} and the $2\pi$-periodicity of $I_n(t)$ and $\theta_n(t)-t$ in $\T$, for all $n\geq 0$ we have
\begin{equation}\label{eq:integrales_cero}
\int_0^{2\pi} p(t)\frac{\partial x}{\partial I}(\theta_n(t),I_n(t))dt=\int_0^{2\pi} p(t)\frac{\partial x}{\partial\theta}(\theta_n(t),I_n(t))dt=0.
\end{equation}
Passing to the limit in the previous equalities we conclude that
\[
\int_0^{2\pi} p(t)\frac{\partial x}{\partial I}(\theta_*+t,I_*)dt=\int_0^{2\pi} p(t)\frac{\partial x}{\partial\theta}(\theta_*+t,I_*)dt=0.
\]
The functions $\frac{\partial x}{\partial I}$ and $\frac{\partial x}{\partial \theta}$ form a fundamental pair of solutions of~\eqref{eq:var}. In consequence $\psi(t,I_*)$ can be expressed as a complex linear combination of them and therefore $\Phi_p(-\theta_*,I_*)=0.$ This is a contradiction with the condition~\eqref{rc}.

Now let us assume that $I_*=+\infty$. In this case, on account of the estimates~\eqref{l111} and~\eqref{EST2} we can ensure that $\theta_n(t)-t\rightarrow \theta_*$ and $I_n(t)\rightarrow+\infty$ uniformly in $t\in\T$. More precisely, $\frac{I_n(t)}{I_n(0)}\rightarrow 1$ uniformly in $t\in\T$. 
On account of Proposition~\ref{prop:limites-uniformes} we have that, as $I\rightarrow+\infty$, 
$\sqrt{I}\frac{\partial x}{\partial I}(\theta,I)\rightarrow \abs{\cos\left(\frac{\theta}{2}\right)}$
uniformly in $\T$ and 
$
\frac{1}{\sqrt{I}}\dot x(\theta,I)\rightarrow-\sin\left(\frac{\theta}{2}\right)\text{sign}\left(\cos\left(\frac{\theta}{2}\right)\right)
$
uniformly in every compact subset of $\T\setminus\{\pi\}$. Thus, on account of $\frac{I_n(t)}{I_n(0)}\rightarrow 1$ uniformly in $t\in\T$, as $n\rightarrow+\infty$,
\[
p(t)\sqrt{I_n(0)}\frac{\partial x}{\partial I}(\theta_n(t),I_n(t))
\rightarrow 
p(t)\abs{\cos\left(\frac{\theta_*+t}{2}\right)}
\]
uniformly in $\T$ and 
\[
p(t)\frac{1}{\sqrt{I_n(0)}}\dot x(\theta_n(t),I_n(t))\rightarrow-p(t)\sin\left(\frac{\theta_*+t}{2}\right)\text{sign}\left(\cos\left(\frac{\theta_*+t}{2}\right)\right)
\]
uniformly in every compact subset of $\T\setminus\{\pi-\theta_*\}$. Consequently, using \eqref{eq:integrales_cero} for this particular case,
\[
0=\int_{-\pi}^{\pi}p(t)\sqrt{I_n(0)}\frac{\partial x}{\partial I}(\theta_n(t),I_n(t))dt
\rightarrow \int_{-\pi}^{\pi} p(t)\abs{\cos\left(\frac{\theta_*+t}{2}\right)}dt
\]
and
\[
0=\int_{-\pi}^{\pi}p(t)\frac{1}{\sqrt{I_n(0)}}\dot x(\theta_n(t),I_n(t))dt\rightarrow-\int_{-\pi}^{\pi}p(t)\sin\left(\frac{\theta_*+t}{2}\right)\text{sign}\left(\cos\left(\frac{\theta_*+t}{2}\right)\right)dt.
\]
In view of Proposition~\ref{prop:cotas-isocrono} we observe that $\frac{1}{\sqrt{I_n(0)}}\dot x(\theta_n(t),I_n(t))$ is uniformly bounded. The convergence of the second integral follows from Lebesgue dominated convergence Theorem. These two integrals vanish and this is not compatible with condition~\eqref{rc}. Indeed we can apply again dominated convergence to prove that
\[
\lim_{r\rightarrow+\infty}\int_0^{2\pi}p(t-\theta)\psi(t,r)dt=\int_0^{2\pi}p(t-\theta)\left(\abs{\cos\left(\frac{t}{2}\right)}+2i\sin\left(\frac{t}{2}\right)\text{sign}\left(\cos\left(\frac{t}{2}\right)\right)\right)dt.
\]
Here we are using~\eqref{eq:expr-var}.
\end{prooftext}

\begin{prooftext}{Proof of Corollary~\ref{corolario}}
The Fourier expansion of $\psi(\cdot,r)$ is of the type
\[
\psi(t,r)=\sum_{m=-\infty}^{+\infty} c_m(r) e^{imt}
\]
and we will derive some properties of the first coefficients. The symmetry $\psi(-t,r)=\overline{\psi(t,r)}$ implies that the numbers $c_0(r)$, $c_1(r)$ and $c_{-1}(r)$ are real. Define
\begin{align*}
d_+(r)&=\frac{1}{2}(c_1(r)+c_{-1}(r))=\frac{1}{2\pi}\int_{-\pi}^{\pi}\left[ Re\psi(t,r)\right]\cos(t) dt,\\
d_-(r)&=\frac{1}{2}(c_1(r)-c_{-1}(r))=\frac{1}{2\pi}\int_{-\pi}^{\pi}\left[ Im\psi(t,r)\right]\sin(t) dt.\\
\end{align*}
From the explicit formula~\eqref{eq:expr-var} we observe that
\begin{equation}\label{d17}
d_-(r)=\frac{1}{2\pi}\int_{-\pi}^{\pi}\frac{\sin^2(t)}{\sqrt{\cos^2(\tfrac{t}{2})+\tfrac{1}{\lambda^4}\sin^2(\tfrac{t}{2})}}dt \geq d_-(0)=\frac{1}{2}.
\end{equation}
Here $\lambda=1+r$. To obtain a lower estimate of $d_+(r)$ we cannot apply the same strategy because the corresponding integrand changes sign. We use the symmetry ($Re\psi(\cdot,r)$ is even) to reduce the interval of integration,
\begin{equation}\label{news}
d_+(r)=\frac{1}{\pi}\int_0^{\pi}\left[Re\psi(t,r)\right]\cos(t)dt=\frac{1}{\pi}\int_0^{\frac{\pi}{2}}\left[Re\psi(t,r)-Re\psi(\pi-t,r)\right]\cos(t)dt.
\end{equation}
A direct computation shows that $Re\dot\psi(t,r)<0$ on $(0,\pi)$. In consequence $Re\psi(\cdot,r)$ is decreasing and $d_+(r)$ is positive. We can go further after observing that $d_+'(r)>0$ for each $r\geq 0$. This is a consequence of long and direct computations coming from~\eqref{eq:expr-var}. Differentiating with respect to $r$,
\begin{equation}\label{magic}
\frac{\partial}{\partial r}\left[Re\psi(t,r)\right]=\frac{\sqrt{2}\bigl(4\sin^4(\tfrac{t}{2})+3\lambda^4\sin^2(\tfrac{t}{2})\bigr)}{\lambda^3\bigl(2+(\lambda^4-1)(1+\cos(t))\bigr)^{\frac{3}{2}}}.
\end{equation}
Next we differentiate this formula with respect to $t$ by quotient rule and observe that all summands on the numerator remain positive on $(0,\pi)$. Then $\frac{\partial}{\partial r}\left[Re\dot\psi(t,r)\right]>0$ and $\frac{\partial}{\partial r}\left[Re\psi(\cdot,r)\right]$ is increasing on $(0,\pi)$. Differentiating with respect to $r$ in~\eqref{news} and looking at the integral over $(0,\tfrac{\pi}{2})$ we conclude that $d_+'(r)<0$. Hence,
\begin{equation}\label{d23}
d_+(r)>d_+(+\infty)=\frac{1}{\pi}\int_0^{\pi}\abs{\cos(\tfrac{t}{2})}\cos(t)dt = \frac{2}{3\pi}.
\end{equation}
The computation of $d_+(+\infty)$ involves a passage to the limit $(r\rightarrow+\infty)$ in the integral. This is justified by the estimate $\abs{Re\psi(t,r)}\leq 2$ implying that there is dominated convergence.

With respect to the coefficient
\[
c_0(r)=\frac{1}{2\pi}\int_{-\pi}^{\pi}Re\psi(t,r)dt
\]
we differentiate in $r$ and apply~\eqref{magic} to deduce that $c_0'(r)>0$. Then, using again dominated convergence,
\[
0\leq c_0(r)<c_0(+\infty)=\frac{1}{\pi}\int_0^{\pi}\abs{\cos(\tfrac{t}{2})}dt=\frac{2}{\pi}.
\]

After these preliminary estimates we are ready to prove that the condition~\eqref{New} implies~\eqref{rc}. From the definition of $\Phi_p(\theta,r)$,
\[
2\pi\Phi_p(\theta,r)=a_0 c_0(r)+d_+(r)(a_1\cos(\theta)-b_1\sin(\theta))+d_-(r)(b_1\cos(\theta)+a_1\sin(\theta))i.
\]
This integral can be identified to a vector in $\R^2$ defined by $w+DR[\theta]v$ where
\[
v=\left(\begin{matrix} a_1\\b_1\\ \end{matrix}\right),\ 
w=\left(\begin{matrix} a_0 c_0(r)\\0\\ \end{matrix}\right),
D=\left(\begin{matrix} d_+(r) & 0 \\0 & d_-(r)\\ \end{matrix}\right),
\]
and $R[\theta]$ is the matrix of a rotation in the counter-clockwise sense. In the next computations we employ the Euclidean norm in $\R^2$, denoted by $|\cdot|$. From~\eqref{d17} and~\eqref{d23} we deduce
\[
\abs{2\pi\Phi_p(\theta,r)}=\abs{w+DR[\theta]v}=\abs{D(D^{-1}w+R[\theta]v)}\geq \frac{2}{3\pi}\abs{D^{-1}w+R[\theta]v}.
\]
In addition, in view of~\eqref{d17} and~\eqref{d23},
\[
\abs{D^{-1}w+R[\theta]v}\geq \abs{v}-\abs{D^{-1}w}=\sqrt{a_1^2+b_1^2}-\frac{c_0(r)\abs{a_0}}{d_+(r)}\geq \sqrt{a_1^2+b_1^2}-3 \abs{a_0} >0.
\]
The last inequality comes from condition~\eqref{New}. Finally,
\[
\abs{\Phi_p(\theta,r)}\geq \frac{1}{3\pi^2}\bigl(\sqrt{a_1^2+b_1^2}-3 \abs{a_0}\bigr)>0,
\]
as we desired.
\end{prooftext}

\subsection{Proof of Theorem~\ref{thm:multiplicativo}}\label{multiplicativa}

According to~\cite{Pinney} the equation
\begin{equation}\label{epin}
\ddot x + x = \frac{\lambda}{x^3},\ x>0
\end{equation}
can be solved explicitly. The solution satisfying $x(0)=x_0>0$, $\dot x(0)=y_0$ is
\[
x(t)=\sqrt{(x_0\cos(t)+y_0\sin(t))^2+\frac{\lambda}{x_0^2}\sin^2(t)}.
\]
For $\lambda$ positive this solution is defined everywhere and it is possible to consider the map
\[
\Phi_{\lambda}:D\rightarrow D,\ (x(0),\dot x(0))\mapsto \left(x\bigl(\frac{\pi}{2}\bigr),\dot x\bigl(\frac{\pi}{2}\bigr)\right)
\]
where $D=(0,+\infty)\times\R$. After some computations we obtain
\[
\Phi_{\lambda}(x_0,y_0)=\left(x_0\varphi_{\lambda}(x_0,y_0),\frac{-y_0}{\varphi_{\lambda}(x_0,y_0)}\right)\ \text{ with }\varphi_{\lambda}(x_0,y_0)=\sqrt{\frac{y_0^2}{x_0^2}+\frac{\lambda}{x_0^4}}.
\]
The solutions of the equation~\eqref{acw} are also defined everywhere and the Poincar\'e map $P:D\rightarrow D$ is given by $P=\Phi_{c}\circ\Phi_1$. Here we are using that the equation~\eqref{epin} is autonomous. Some additional computations show that
\[
P(x_0,y_0)=\left( x_0\Pi(x_0,y_0), \frac{y_0}{\Pi(x_0,y_0)}\right)
\]
where $\Pi(x_0,y_0)=\sqrt{\frac{x_0^2y_0^2+c}{x_0^2y_0^2+1}}$. Due to the structure of $P$ it is easy to check that the function $I(x_0,y_0)=x_0y_0$ is a first integral, meaning that $I\circ P = I$. In consequence $\Pi$ is also a first integral and we can compute explicitly the orbits of $P$, namely $(x_n,y_n)=P^n(x_0,y_0)=\left(x_0\Pi(x_0,y_0)^n,\frac{y_0}{\Pi(x_0,y_0)^n}\right)$. When $c>1$ the function $\Pi$ takes values above $1$ and so the sequence $x(n\pi)=x_n$ is a geometric progression going to infinity. When $c<1$ this is also the case of $\dot x(n\pi)=y_n$. In both cases all solutions are unbounded.
\Qed

\section{Appendix: a result on isochronous centers}

Next result is concerned with the behaviour at infinity of isochronous center when it presents a single asymptote and it is based in \cite[Lemma 8]{BFS}. Here we present a more general version by considering the convexity located at the endpoints of the interval of definition of $V$ instead of the global convexity.

\begin{lema}\label{lema:asymptota}
Let $V:(a,+\infty)\rightarrow[0,+\infty)$ be an analytic $2\pi$-isochronous potential such that $V(0)=V'(0)=0$ and $xV'(x)>0$ if $x\neq 0$. Assume that $\lim_{x\rightarrow a} V(x)=+\infty$ and $\lim_{x\rightarrow +\infty} V(x)=+\infty$. In addition we assume that there exist $\zeta,M>0$ such that $V$ is convex in $(a,a+\zeta)$ with $a+\zeta<0$ and in $(M,+\infty)$. Then $\lim_{x\rightarrow+\infty} \bigl(V'(x)-\frac{x}{4}\bigr) = -a/4$.
\end{lema}

\begin{demo}
Let us denote by $V_-$ and $V_+$ the restriction of $V$ on the negative and positive axis, respectively. From \cite[Theorem B]{CMV} the fact that $V$ is an analytic $2\pi$-isochronous implies that
\begin{equation}\label{eq:iso}
V_+(x)=(x-\sigma(x))^2/8
\end{equation}
for all $x>0$, where $\sigma(x)\!:=V_{-}^{-1}(V_+(x))$. Hence, for $x$ positive, $V_+'(x)=\frac{1}{4}(x-\sigma(x))(1-\sigma'(x))$. Due to $V_+(x)\rightarrow+\infty$ as $x\rightarrow+\infty$ and $V_{-}^{-1}(y)\rightarrow a$ as $y\rightarrow+\infty$ we have that $\sigma(x)\rightarrow a$ as $x\rightarrow+\infty$. On the other hand, $V_-^{-1}$ is a convex decreasing function in $(V_-^{-1}(a+\zeta),+\infty)$ with limit $a<0$ at infinity so by Lemma~\ref{lema:convex} we have that 
\begin{equation}\label{eq:limit}
\lim_{y\rightarrow+\infty}y(V_{-}^{-1})'(y)=0.
\end{equation}

We claim at this point that $\frac{xV_+'(x)}{V_+(x)}$ is bounded as $x\rightarrow+\infty$. Indeed, $V_+$ is a positive increasing convex function in $(M,+\infty)$ so for all $x,y\in(M,+\infty)$ we have $V_+(y)\geq V_+(x)+V_+'(x)(y-x)$. Taking $x>M$ and $y=2x$ we have then
\[
0 \leq \frac{xV_+'(x)}{V_+(x)}\leq \frac{V_+(2x)}{V_+(x)}-1
\]
for all $x>M$.
According with expression in \eqref{eq:iso} and using $\sigma(x)\rightarrow a$ as $x\rightarrow+\infty$, we have
\[
\lim_{x\rightarrow+\infty}\frac{V_+(2x)}{V_+(x)}-1=\lim_{x\rightarrow+\infty}\frac{(2x-\sigma(2x))^2}{(x-\sigma(x))^2}-1=3
\]
which proves the claim.

By the expression of $\sigma(x)$ we have that 
\[
x\sigma'(x)=x(V_-^{-1})'(V_+(x))V_+'(x)
=\frac{xV_+'(x)}{V_+(x)}(V_-^{-1})'(V_+(x))V_+(x).
\]
On account of equality \eqref{eq:limit} and using $V_+(x)\rightarrow+\infty$ as $x\rightarrow+\infty$ we have $\lim_{x\rightarrow+\infty}(V_-^{-1})'(V_+(x))V_+(x)=0$. On the other hand, $\frac{xV_+'(x)}{V_+(x)}$ is bounded as $x\rightarrow+\infty$. Thus, $x\sigma'(x)\rightarrow 0$ as $x\rightarrow+\infty$ (and, particularly, $\sigma'(x)\rightarrow 0$ as $x\rightarrow+\infty$). Consequently,
\[
\lim_{x\rightarrow+\infty}\left(V_+'(x)-\frac{x}{4}\right)=\lim_{x\rightarrow+\infty} \left(-\frac{1}{4}\sigma(x)+\frac{1}{4}\sigma(x)\sigma'(x)-\frac{1}{4}x\sigma'(x)\right)=-\frac{a}{4}
\]
and so the result is proved.
\end{demo}


\begin{thebibliography}{9}

\bibitem{ACW}
J.~Ai, K.S.~Chou and J.~Wei,
\emph{Self-similar solutions for the anisotropic affine curve shortening problem},
Calc. Var. Partial Differential Equations 13 (2001) 311--337.

\bibitem{AO}
J.M.~Alonso and R.~Ortega,
\emph{Roots of unity and unbounded motions of an asymmetric oscillator},
J. Differential Equations 143 (1998) 201--220.

\bibitem{BonFab}
D.~Bonheure and C.~Fabry,
\emph{Littlewood's problem for isochronous oscillators},
Arch. Math. (Basel) 93 (2009) 379--388.

\bibitem{BFS}
D.~Bonheure, C.~Fabry and D.~Smets,
\emph{Periodic solutions of forced isochronous oscillators at resonance},
Discrete Contin. Dyn. Syst. 8 (2002) 907--930.

\bibitem{BES}
B.M.~Brown, M.S.P.~Eastham and K.M.~Schmidt,
\emph{Periodic differential operators}, Operator Theory: Advances and Applications 230, Birkh\"auser/Springer Basel AG, Basel (2013).

\bibitem{CapDamLiu}
A.~Capietto, W.~Dambrosio and B.~Liu, 
\emph{On the boundedness of solutions to a nonlinear singular oscillator},  
Z. Angew. Math. Phys. 60 (2009) 1007--1034.

\bibitem{CMV}
A.~Cima, F.~Ma\~nosas and J.~Villadelprat,
\emph{Isochronicity for several classes of Hamiltonian systems},
J. Differential Equations 157 (1999) 373--413.

\bibitem{Dancer}
E.N.~Dancer,
\emph{Boundary value problems for weakly nonlinear ordinary differential equations},
Bull. Aust. Math. Soc. 15 (1976) 321--328.

\bibitem{Kamke}
E.~Kamke,
\emph{A new proof of Sturm's comparison theorems},
Amer. Math. Monthly 46 (1939) 417--421.

\bibitem{LiJin}
X.~Li and S.~Jin, 
\emph{Boundedness in forced isochronous oscillators}, 
J. Math. Anal. Appl. 460 (2018) 714--736.

\bibitem{Liu}
B.~Liu,
\emph{Boundedness in asymmetric oscillators},
J. Math. An. App. 231 (1999) 355--373.

\bibitem{NewLiu}
B.~Liu, 
\emph{Quasi-periodic solutions of forced isochronous oscillators at resonance}, 
J. Differential Equations 246 (2009) 3471--3495.

\bibitem{Lloyd}
N.G.~Lloyd,
\emph{Degree theory}, 
Cambridge Tracts in Mathematics 73. Cambridge University Press, Cambridge-New York-Melbourne (1978).

\bibitem{ManRojVil2015}
F.~Ma{\~n}osas, D.~Rojas and J.~Villadelprat,
\emph{The criticality of centers of potential systems at the outer boundary}, J. Differential Equations 260 (2016) 4918--4972.

\bibitem{Massera}
J.~Massera,
\emph{The existence of periodic solutions of systems of differential equations},
Duke Math. J. 17 (1950), 457--475.

\bibitem{Ortega96}
R.~Ortega,
\emph{Asymmetric oscillators and twist mappings}, J. London Math. Soc. 53 (1996) 325--342.

\bibitem{Ortega}
R.~Ortega,
\emph{Periodic perturbations of an isochronous center},
Qualitative Theory of Dynamical Systems 3 (2002) 83--91.

\bibitem{Ortega3}
R.~Ortega,
\emph{Periodic differential equations in the plane: a topological perspective}, to appear.

\bibitem{OrtRoj}
R.~Ortega and D.~Rojas,
\emph{A proof of Bertrand's theorem using the theory of isochronous potentials}, to appear in J. Dyn. Diff. Equat. (2018).

\bibitem{Pinney}
E.~Pinney,
\emph{The nonlinear differential equation $y''+p(x)y+cy^{-3}=0$},
Proc. Amer. Math. Soc. 1 (1950) 681.

\bibitem{Urabe62}
M.~Urabe,
\emph{The potential force yielding a periodic motion whose period is an arbitrary continuous function of the amplitude of the velocity},
M. Arch. Rational Mech. Anal. 11 (1962) 27--33.

\end{thebibliography}
\end{document}